\documentclass[12pt]{article}

\usepackage{amssymb,amsmath,amsthm,hyperref,amssymb,latexsym,color,titlesec,enumerate}

\usepackage[shortlabels]{enumitem}
\usepackage{rotating} 

\usepackage{authblk}

\theoremstyle{plain}
\newtheorem{theorem}{Theorem}[section]
\newtheorem{proposition}[theorem]{Proposition}
\newtheorem{lemma}[theorem]{Lemma}
\newtheorem{corollary}[theorem]{Corollary}

\theoremstyle{definition}
\newtheorem{definition}[theorem]{Definition}
\newtheorem{remark}[theorem]{Remark}
\newtheorem{examplex}[theorem]{Example}
	\newenvironment{example}
	{\pushQED{\qed}\examplex}
	{\popQED\endexamplex}

\def\g{\mathfrak{g}}
\def\G{\mathcal{G}}
\DeclareMathOperator{\dist}{dist}

\newcommand{\be}{\begin{equation}}

\newcommand{\ee}{\end{equation}}

\titleformat{\subsubsection}
{\normalfont\normalsize\bfseries}{\thesubsubsection}{1em}{}
\titlespacing*{\subsubsection}
{0pt}{0.5ex}{-0.5ex plus 0ex}

\begin{document}

\title{\bf{Rotationally Symmetric Extremal K\"ahler Metrics on $\mathbb C^n$ and $\mathbb C^2\setminus \{0\}$ }} 
\renewcommand\footnotemark{}

\author{Selin Ta\c skent}

\date{}
\maketitle

\begin{abstract}
	In this paper, we study rotationally symmetric extremal K\"ahler metrics on $\mathbb C^n$ ($n\geq 2$) and $\mathbb C^2 \backslash\{0\}$. We present a classification of such metrics based on the zeros of the polynomial appearing in Calabi's Extremal Equation. As applications, we prove that there are no $U(n)$ invariant complete extremal K\"ahler metrics on $\mathbb C^n$ with positive bisectional curvature, and we give a smooth extension lemma for $U(n)$ invariant extremal K\"ahler metrics on $\mathbb{C}^n\backslash\{0\}$. We retrieve known examples of smooth or singular extremal K\"ahler metrics on Hirzebruch surfaces, bundles over $\mathbb{CP}^1$, and weighted complex projective spaces. We also show that certain solutions on $\mathbb C^2\backslash\{0\}$ correspond to new complete families of constant-scalar-curvature K\"ahler  and strictly extremal K\"ahler metrics on complex line bundles over $\mathbb{CP}^{1}$ and on $\mathbb C^2\backslash\{0\}$. 
\end{abstract}


\pagenumbering{arabic}

\section{Introduction}

In this paper, we study rotation invariant K\"ahler metrics on $\mathbb C^n$ ($n\geq 2$) and $\mathbb C^2 \backslash\{0\}$ that are solutions to Calabi's extremal equation  \cite{Calabi82}. The rotational symmetry allowed Calabi to reduce the extremal equation to an ordinary differential equation for the K\"ahler potential, $u$, written as a function of $s=|z|^2$.  In fact,
\begin{equation} \label{ODE}
sg'(s)=F(g(s)) \textrm{ where } g(s)=su'(s).
\end{equation} 


Calabi defined $\mathcal{F}_k^n$ to be the $k$-twisted $\mathbb{CP}^1$-bundle over $\mathbb{CP}^{n-1}$ and  defined a biholomorphism 
\be
\tilde{p} : (\mathbb{C}^n \backslash \{0\})/\mathbb{Z}_k \rightarrow \mathring{\mathcal{F}}_k^n,
\ee
where $\mathring{\mathcal{F}}_k^n$ is the complement of $|z|=0$ and $|z|=\infty$ sections.
He derived conditions on $F$ of (\ref{ODE}) for adding  $\mathbb{CP}^{n-1}$ smoothly to $(\mathbb{C}^n \backslash \{0\})/\mathbb{Z}_k$ at $|z|=0$ and $|z|=\infty$.   Cao \cite{C94} used the map $\tilde{p}$ to produce $U(n)$ invariant, complete gradient K\"ahler-Ricci soliton metrics (GKRS) on line bundles over  $\mathbb{CP}^{n-1}$. Feldman-Ilmanen-Knopf \cite{FIK03} generalized this approach by producing $U(n)$-invariant GKRS metrics on $(\mathbb{C}^n \backslash \{0\})/\mathbb{Z}_k$, where they allowed new boundary behavior at $z=0$ and $|z|=\infty$, such as adding an orbifold point or a singular $\mathbb{CP}^{n-1}$. In this method, completeness of the metric as $|z|\rightarrow 0$ (resp. $|z|\rightarrow \infty$) can be checked by showing that the integral that gives geodesic distance to $|z|=0$ (resp. $|z|=\infty$)  diverges. (See Section \ref{section2}).

In \cite{HeLi}, He and Li present a classification of solutions to Equation \eqref{ODE}, that are constant scalar curvature K\"ahler (cscK) metrics on $\mathbb C^n$, $\mathbb C^2 \backslash\{0\}$, and $\mathbb C^3 \backslash\{0\}$  based on the zeros of the function $F$ (see Theorems \ref{n smooth solution}, \ref{2 zero singular solution}, \ref{2 negative singular solution}). In this paper, we extend their results to give a classification of solutions to the extremal equation \eqref{ODE} on  $\mathbb C^n$ and $\mathbb C^2 \backslash\{0\}$ (see Theorems \ref{Thm2.2}, \ref{lemma2}). This classification contains only those equations that give rise to K\"ahler potentials defined on the interval $(0,\infty)$, hence to extremal K\"ahler metrics defined on $\mathbb C^n$ and $\mathbb C^2 \backslash\{0\}$. Existence of K\"ahler potentials that appear in our list is guaranteed by technical lemmas given by He-Li (see Lemmas \ref{lemma6.2}, \ref{lemma6.1}).
 
In Theorem \ref{Thm2.2}, we classify solutions of Equation \eqref{ODE} that give extremal K\"ahler metrics on $\mathbb{C}^n$, based on zeros of $F$. Though other classifications were obtained by mathematicians (e.g. Bryant \cite{B01}, Gauduchon \cite{G19}), we include it here for comparison and for its applications (see Corollary \ref{thm1}, Lemma \ref{addingsmoothpt}, Proposition \ref{prop2.8}, and Remark \ref{rmk-BKL}).

 \begin{theorem}\label{Thm2.2}
 	Let $n\geq 2$ and $u:[0,\infty)\rightarrow \mathbb R$ be 
 	the potential of an extremal K\"ahler metric 
 	on $\mathbb C^n$. Then, for $g(s) = su'(s)$, one of the following is true:
 	\begin{enumerate}
 		\item \label{cscK}$\omega$ is a cscK metric.
 		\item \label{metric beta}There exist constants $\beta,c$ with $\beta>0$ such that $g(s)$ is the smooth strictly increasing function $g:(0,\infty) \rightarrow(0,\beta)$ determined by the ordinary differential equation
 		$$sg'(s)=F(g(s)) = \frac{1}{\beta^2}g(s) (g(s)-\beta)^2.$$ 
 		These metrics have positive non-constant scalar curvature.
 		
 		\item \label{metric betagamma}There exist constants $\gamma,\beta,c$ with $\gamma<0<\beta$ such that  $g(s)$ is the smooth strictly increasing function ranging from $g:(0,\infty)\rightarrow(0, \beta)$ determined by the ordinary differential equation
 		$$sg'(s)=F(g(s)) = \frac{1}{\beta\gamma}g(s) (g(s)-\beta)(g(s)-\gamma).$$

 		\item \label{metric betagamma2} There exist constants $\gamma,\beta,c$ with $0<\beta<\gamma$ such that  $g(s)$ is the smooth strictly increasing function  $g:(0,\infty)\rightarrow(0, \beta)$ determined by the ordinary differential equation
 		$$sg'(s)=F(g(s)) = \frac{1}{\beta\gamma}g(s) (g(s)-\beta)(g(s)-\gamma).$$
 		
 	\end{enumerate}
 \end{theorem}

In Theorem \ref{lemma2}, we give a classification of solutions of Equation \eqref{ODE} with non-constant scalar curvature on $\mathbb{C}^2\backslash\{0\}$. Theorem \ref{lemma2}, together with scalar-flat and positive constant scalar curvature K\"ahler metrics given by He-Li (Theorems \ref{2 zero singular solution}, \ref{2 negative singular solution}), gives the complete list of solutions to Calabi's extremal equation \eqref{ODE} on $\mathbb{C}^2\backslash\{0\}$.

Naturally, our classifications (Theorems \ref{Thm2.2}, \ref{lemma2}) include examples that have been constructed before by mathematicians seeking  scalar-flat, cscK, self-dual K\"ahler ($n=2$), Bochner-K\"ahler, and extremal K\"ahler metrics using other methods. In fact, Theorem \ref{Thm2.2} provides a classification where each case corresponds to a known extremal K\"ahler metric. These metrics live on  $\mathring{\mathcal{F}}_k^n$ via the map $\tilde{p}$, and are completed in accordance with their boundary behavior as $|z|\rightarrow 0$ and $|z|\rightarrow \infty$. (See Remark \ref{rmk-BKL}). Theorem \ref{lemma2} gives Calabi's extremal K\"ahler metrics on Hirzebruch surfaces $\mathcal{F}_k^2$, $k\geq 1$; positive cscK metrics on $\mathcal{O}_{\mathbb{CP}^1}(k)$, $k\geq 1$, (Example \ref{ex4.6}); strictly extremal K\"ahler metrics on $\mathcal{O}_{\mathbb{CP}^1}(-k)$, $k\geq 1$, (Example \ref{example 3.7}); and strictly extremal K\"ahler metrics on $\mathbb{C}^2\backslash \{0\}$ (Example \ref{ex4.8}). To the best of our knowledge, the last three metrics are new.

In this paper, we pursue the methods used in \cite{Calabi82,HeLi, C94, FIK03} to produce new and well-known examples of scalar-flat, cscK, and strictly extremal metrics. We use Theorem \ref{Thm2.2} to obtain a \textit{smooth extension lemma} on $\mathbb{C}^n \backslash\{0\}$  (Lemma \ref{addingsmoothpt}) and to obtain a result on non-positivity of bisectional curvature (Proposition \ref{prop2.8}).  

\subsubsection*{Background}

In Riemannian geometry, the curvature tensor can be decomposed into its irreducible components under the action of the orthogonal group. The three components are the scalar curvature, the traceless Ricci curvature and the Weyl curvature. In K\"ahler geometry, there is a similar decomposition under the action of the unitary group. The K\"ahler decomposition has the scalar curvature, the traceless Ricci curvature and the Bochner curvature as its components.  For a K\"ahler manifold of real dimension four ($n=2$), the anti-self-dual part of the Weyl tensor is the same as the Bochner tensor. Hence, in real dimension four, Bochner-K\"ahler metrics and self-dual K\"ahler metrics are the same. Bochner-K\"ahler metrics are known to be extremal in the sense of Calabi. It is easy to see that $U(n)$ invariant extremal K\"ahler metrics on $\mathbb{C}^n\backslash\{0\}$ extend smoothly to the origin if and only if they are Bochner-K\"ahler. (we include this as an equivalent condition in our smooth extension lemma (Lemma \ref{addingsmoothpt})). The first examples of $U(n)$ invariant Bochner-K\"ahler metrics with non-constant scalar curvature were found by Tachibana and Liu \cite{TL70}, using a K\"ahler potential defined on an interval $s\in I\subset \mathbb R$, ($s=|z|^2$).

In \cite{TL70}, Tachibana-Liu  give the ansatz for $U(n)$ invariant  Bochner-K\"ahler metrics on $\mathbb{C}^n$.  In \cite{B01}, Bryant gives a classification of these metrics, and determines that $U(n)$ invariant complete Bochner-K\"ahler metrics with non-constant scalar curvature are unique (up to scaling). Furthermore, he gives a complete classification of Bochner-K\"ahler metrics in dimension $n$. He shows that when singularities are allowed, every weighted projective space supports a Bochner-K\"ahler metric. 

In \cite{AG02}, Apostolov-Gauduchon obtain a local classification of self-dual Hermitian-Einstein metrics, and this implies a local classification of self-dual K\"ahler metrics. They prove that such metrics have local cohomogeneity at most $2$. Their classification is an independent alternative to Bryant's classification when $n=2$.

In \cite{G19}, Gauduchon uses  momentum profile  to study Calabi's $U(n)$ invariant extremal K\"ahler metrics on $\mathbb{C}^n \backslash \{0\}$. We refer the reader to Hwang-Singer \cite{HS02} for an introduction on momentum profile approach, which explains that it is possible to solve the equations explicitly after a change in coordinates. Gauduchon  gives two examples of $U(n)$ invariant extremal K\"ahler metrics on $\mathbb{C}^n$ with non-constant scalar curvature.  The first metric is due to Gauduchon \cite{G19}, where $\mathbb{C}^n$ is completed by adding a singular $\mathbb{CP}^{n-1}$ at $|z|=\infty$ with cone angles $2\pi\theta$, $0<\theta<1$. (See Example \ref{ex gaud} for details). The second metric is complete and was first discovered by Olivier Biquard. This metric is discussed in more detail in Example \ref{biquard}. These two metrics appear as Cases \eqref{metric beta} and \eqref{metric betagamma2} of Theorem \ref{Thm2.2}. Moreover, Case \eqref{metric betagamma} of Theorem \ref{Thm2.2} gives incomplete extremal (Bochner-) K\"ahler metrics on $\mathbb{C}^n$. In Example \ref{ex gkn}, we show that these metrics induce extremal K\"ahler \textit{orbifold} metrics on $G_k^n = \mathbb{CP}^n_{[k,1,\dots,1]}$. These orbifold metrics were studied by other mathematicians in various settings (see Remark \ref{rmk-BKL}). We note that together with Euclidean and Fubini-Study metrics given by Case \eqref{cscK} of Theorem \ref{Thm2.2}, these three metrics exhaust our list on $\mathbb{C}^n$.

The orbifold metrics on $G_k^2 = \mathbb{CP}^2_{[k,1,1]}$ can be found implicitly in the earlier work \cite{HS97} of Hwang and Simanca. In their paper, they show that extremal K\"ahler metrics on Hirzebruch surfaces which are locally conformally equivalent to Einstein metrics are exactly the ones that are degenerate, \textit{i.e.} the ones that would give smooth orbifold metrics on $G_k^2$ after blowdown of zero section in $\mathcal{F}_k^2$, $k\geq 1$. Furthermore, for the extremal representative of each K\"ahler class in $\mathcal{F}_k^2$, they give conditions for the scalar curvature to be everywhere positive, or to vanish on specific subsets of $\mathcal{F}_k^2$. Taking the complement of the zero set of scalar curvature in $\mathcal{F}_k^2$, Hwang and Simanca give extremal K\"ahler metrics on disk bundles over $\mathbb{CP}^1$.

In \cite{LeBrun88}, LeBrun constructed scalar-flat ALE  K\"ahler metrics on 
$\mathcal{O}_{\mathbb{CP}^{1}}(-k)$ for $k\geq 1$ as counterexamples to \textit{Positive Action Conjecture}. For $k=1$, this is the Burns metric \cite{LeBrun88}. For $k=2$, this metric is the Ricci-flat Eguchi-Hanson metric \cite{EH79}. (See Remark \ref{remark 3.4}). Simanca \cite{S91} constucted scalar-flat K\"ahler metrics on $\mathcal{O}_{\mathbb{CP}^{n-1}}(-1)$ ($n>2$). Scalar-flat metrics on $\mathcal{O}_{\mathbb{CP}^{n-1}}(-k)$ ($k,n>2$) are given by Pedersen-Poon \cite{PP91}. In \cite{AR17}, Apostolov and Rollin generalize the construction of scalar-flat ALE metrics on $\mathcal{O}_{\mathbb{CP}^{n-1}}(-k)$  to non-compact weighted projective spaces $\mathbb{CP}^{n}_{[-a_0,a_1,\dots,a_n]}$. In the special case of  $\mathbb{CP}^{n}_{[-k,1,\dots,1]}$, their metrics coincide with the aforementioned scalar-flat metrics. In \cite{G19}, Gauduchon uses momentum profile to study $U(n)$ invariant scalar-flat K\"ahler metrics on $\mathbb{C}^n\backslash\{0\}$, and provides a method for obtaining the metrics mentioned above. In Remark \ref{remark 3.4} we give the list of scalar-flat metrics induced by Theorem \ref{2 zero singular solution} (Theorem 1.2 in \cite{HeLi}) on $\mathcal{O}_{\mathbb{CP}^1}(-k)$, $k\geq 1$.

$U(2)$ invariant extremal K\"ahler metrics with positive constant scalar curvature on $\mathbb{C}^2\backslash\{0\}$ are listed by He and Li (see Theorem \ref{2 negative singular solution}). For a list of smooth and singular metrics induced by Theorem \ref{2 negative singular solution} on $\mathbb{CP}^2$, Hirzebruch surfaces, and $\mathcal{O}_{\mathbb{CP}^1}(k)$, $k\geq 1$,  see Remark \ref{rmk-pos-csck}. To the best of our knowledge, metrics with positive constant-scalar-curvature on $\mathcal{O}_{\mathbb{CP}^1}(k)$, $k\geq 1$, do not appear in the literature.

The extremal K\"ahler metrics induced by Theorem \ref{lemma2} include Calabi's extremal metrics with non-constant scalar curvature on Hirzebruch surfaces, complete strictly extremal metrics on $\mathcal{O}_{\mathbb{CP}^1}(-k)$, $k\geq 1$, and $\mathbb{C}^2\backslash\{0\}$. The last two metrics seem to be new.

Since our classification only includes solutions defined on all of $\mathbb C^2 \backslash\{0\}$, it does not include examples of scalar-flat, cscK, and strictly extremal metrics on \textit{disk bundles} over $\mathbb{CP}^1$ given by Simanca, Hwang-Simanca, and Abreu in \cite{S91, HS97, Abreu10}.   

\subsubsection*{Results}

In \cite{HeLi}, He and Li give the list of $U(n)$ invariant cscK metrics on $\mathbb{C}^n$ (See Theorem \ref{n smooth solution}). In Theorem \ref{Thm2.2}, we extend their approach to extremal K\"ahler metrics with non-constant scalar curvature on  $\mathbb{C}^n$. Theorems \ref{n smooth solution} and \ref{Thm2.2} together give a complete list of solutions to extremal equation $sg'(s)=F(g(s))$ on $(0,\infty)$, where $F$ is given as in \eqref{odesmooth}. 

As a straightforward application, we prove that there do not exist $U(n)$ invariant extremal K\"ahler metrics with negative scalar curvature on $\mathbb{C}^n$ (Corollary \ref{thm1}). We show that a $U(n)$ invariant extremal K\"ahler metric on $\mathbb C^n \backslash\{0\}$ can be extended smoothly to the origin if and only if it satisfies one of the equivalent conditions given in Lemma \ref{addingsmoothpt}. This lemma improves existing results by reducing the smoothness condition to a single equation, namely $F(0)=0$ (equivalently, $\displaystyle\lim_{s\rightarrow 0^+} g(s)=0$). In Proposition \ref{prop2.8}, we use Wu and Zheng's characterization of $U(n)$ invariant K\"ahler metrics  on $\mathbb{C}^n$ with positive bisectional curvature to show that none of these metrics can be extremal in the sense of Calabi.

In Section \ref{section examples}, we retrieve known complete extremal K\"ahler metrics on $\mathbb{C}^n$, $\mathbb{CP}^n$, and $G_k^n=\mathbb{CP}_{[k,1,\dots,1]}^n$, that are induced by solutions given in Theorem \ref{Thm2.2}.

In \cite{HeLi}, He and Li give the list of scalar-flat and positive constant-scalar-curvature K\"ahler metrics defined on $\mathbb C^2 \backslash\{0\}$.  Here we extend their approach to \emph{extremal} K\"ahler metrics with non-constant scalar curvature on $\mathbb C^2 \backslash\{0\}$. Together with He and Li's Theorems \ref{2 zero singular solution}, \ref{2 negative singular solution}, our Theorem \ref{lemma2} classify solutions on $\mathbb{C}^2\backslash\{0\}$ to the extremal K\"ahler equation $sg'(s)=F(g(s))$, based on zeros of $F$.

The family of $U(2)$ invariant extremal K\"ahler metrics on $\mathbb{C}^2\backslash\{0\}$ given by Theorems \ref{2 zero singular solution}, \ref{2 negative singular solution}, and \ref{lemma2} can be used to write down  extremal K\"ahler metrics on line bundles over $\mathbb{CP}^1$. 

In Remark \ref{remark 3.4}, we retrieve well known examples of scalar flat metrics on $\mathcal{O}_{\mathbb{CP}^1}(-k)$, $k\geq 1$, induced by solutions given in Theorem \ref{2 zero singular solution}. Similarly, in Remark \ref{rmk-pos-csck}, we obtain examples of positive cscK metrics on $\mathbb{CP}^2$, Hirzebruch surfaces $\mathcal{F}_k^2$ (with cone singularities), and on $\mathcal{O}_{\mathbb{CP}^1}(k)$, $k\geq 1$, induced by the solutions given in Theorem \ref{2 negative singular solution}. To the best of our knowledge, positive cscK metrics on $\mathcal{O}_{\mathbb{CP}^1}(k)$, $k\geq 1$, do not appear in the literature. (See Example \ref{ex4.6}).

Some of the metrics induced by the solutions listed in Theorem \ref{lemma2} include Calabi's extremal metrics with non-constant scalar curvature on Hirzebruch surfaces, complete strictly extremal metrics on $\mathcal{O}_{\mathbb{CP}^1}(-k)$, $k\geq 1$, (see Example \ref{example 3.7}), and $\mathbb{C}^2\backslash\{0\}$ (see Example \ref{ex4.8}). The last two metrics seem to be new.

\subsubsection*{Acknowledgments}

I would like to thank my doctoral advisor Xiu-Xiong Chen for his guidance and for suggesting that I extend the classification of He and Li to extremal K\"ahler metrics, give an improved smooth extension lemma, and utilize closing conditions to complete these metrics to obtain new examples.   I would like to thank Claude LeBrun for sharing his invaluable insight on the subject during preparation of my thesis and for mentioning the works of other mathematicians. I'd like to thank Christina Sormani for encouraging me to publish this paper.

\section{$U(n)$ invariant K\"ahler metrics on $\mathbb{C}^n\backslash\{0\}$ } \label{section2}

In this section, we present a summary of methods introduced by Calabi \cite{Calabi82}, Cao \cite{C94}, Wu-Zheng \cite{WZ11}, Feldman-Ilmanen-Knopf \cite{FIK03}, and He-Li \cite{HeLi}, which we use throughout the paper.

\subsection{K\"ahler Potentials}\label{section2.1}
We begin with a well-known decomposition of $U(n)$ invariant K\"ahler metrics on $\mathbb{C}^n\backslash\{0\}$, as introduced in \cite{HeLi} and \cite{FIK03}.

Let $u(s):(0,\infty)\rightarrow \mathbb R$ be a smooth function where $s=|z|^2$. Then, the real $(1,1)$-form
\begin{equation}\label{metric}
	\omega = i\partial\overline\partial u =
	i\sum_{j,k=1}^{n} (\delta_{jk} u'(s)+u''(s)\overline z_j z_k ) dz^j \wedge d\overline z^k
\end{equation} 
gives a positive definite K\"ahler metric on $\mathbb{C}^n\backslash \{0\}$ if and only if
\begin{equation}\label{pos.eqn.}
	u'(s)>0,\qquad u'(s)+su''(s)>0.
\end{equation}

We introduce the function $g(s)=su'(s)$ and reformulate \eqref{pos.eqn.} as
\begin{equation}\label{pos.eqn.g}
	g(s)>0,\qquad g'(s)>0.
\end{equation}
We note that the function $g:(0,\infty)\rightarrow\mathbb{R}$ satisfying \eqref{pos.eqn.g} is positive and strictly increasing. Therefore, $\displaystyle\lim_{s\rightarrow 0^+} g(s)=A$ and $\displaystyle\lim_{s\rightarrow \infty} g(s)=B$ always make sense. We also see that $0\leq A <B \leq +\infty$.

Let us write the metric \eqref{metric} in the form
\begin{align}\label{M1}
\nonumber	\g &= \left(\frac{1}{s} g(s) \delta_{jk} + \frac{1}{s^2}(sg'(s)-g(s)) \overline{z}_j z_k\right) dz^j \otimes d\overline{z}^k\\
&=\left( g(s) \left(\frac{1}{s}\delta_{jk} - \frac{1}{s^2}\overline{z}_j z_k \right)  +sg'(s) \left(\frac{1}{s^2} \overline{z}_j z_k\right)  \right)  dz^j \otimes d\overline{z}^k
\end{align}

We will view $\mathbb{CP}^{n-1}$ as the quotient $(\mathbb{C}^n\backslash\{0\})\slash \mathbb{C}^*$ as well as the quotient $S^{2n-1}(1)\slash S^1$. Let $\pi_1$ and $\pi_2$ denote the corresponding projection maps onto $\mathbb{CP}^{n-1}$, respectively.

The real part of the standard Hermitian product on $\mathbb{C}^n$ induces the Riemannian metric on $S^{2n-1}(1)$. The standard Fubini-Study metric $\g_{FS}$ on $\mathbb{CP}^{n-1}$ is induced by the Riemannian submersion $\pi_2: S^{2n-1}(1)\rightarrow \mathbb{CP}^{n-1}$.

The metric \eqref{M1} can be expressed as in \cite{FIK03} by
\begin{align}\label{M2}
\g &=g(s)(\g_{S^{2n-1}} -\eta\otimes\eta) + sg'(s) \left(\frac{1}{4s^2}ds\otimes ds + \eta\otimes\eta\right)\\
\nonumber&=g(s) \pi_1^* \g_{FS} + sg'(s) \g_{cyl}.
\end{align}
Here $\eta$ gives the $1$-form $d\theta$ when restricted to each complex line through the origin.

Let us introduce the new parameter $r=\sqrt{s}$ and write $$sg'(s)\g_{cyl} = g'(s)(dr\otimes dr + r^2 d\theta\otimes d \theta)$$ on a complex line through the origin. Note that straight lines through the origin coincide with minimal geodesics of the $U(n)$-invariant metric $\g$. It follows that geodesic distance from $z=0$ to $z$ is given by
\begin{equation}
\tilde{r} = \dist(0,z) = \frac{1}{2} \int_0^{s} \sqrt{\frac{g'(s)}{s}}\ ds.
\end{equation}
where $s=|z|^2$. We note that $g'(s) dr\otimes dr = d\tilde{r}\otimes d\tilde{r}$.

We also note that a metric $\g$ on $\mathbb{C}^n$ given by \eqref{M1} is complete if and only if $\displaystyle\int_0^\infty \sqrt{\frac{g'(s)}{s}}\ ds = \infty.$

\subsection{Positive Bisectional Curvature}\label{section2.2}
In \cite{K77}, Klembeck computed the components of the curvature tensor with respect to the orthonormal frame $\{e_1 = \frac{1}{\sqrt{g'}}\ \partial z_1 , e_2=\frac{1}{\sqrt{u'(s)}}\ \partial z_2 , \dots, e_n=\frac{1}{\sqrt{u'(s)}}\ \partial z_n\} $
at a fixed point $(z_1,0,\dots,0)$. 

The nonzero terms are denoted by $A,B,C$ and are given as follows. ($2\leq i\neq j\leq n$)
\begin{align*}
A &=R_{1\overline{1}1\overline{1}} = -\frac{1}{g'} \left(\frac{sg''}{g'}\right)'\\
B &=R_{1\overline{i}1\overline{i}} = \frac{u''}{(u')^2} - \frac{g''}{u'g'}\\
C &=R_{i\overline{i}i\overline{i}} = 2 R_{i\overline{i}j\overline{j}}  = -\frac{2u''}{(u')^2}
\end{align*}

\begin{theorem}[Wu-Zheng \cite{WZ11}] Let $\g$ be a complete $U(n)$ invariant K\"ahler metric on $\mathbb{C}^n$ ($n\geq 2$). Then $\g$ has positive bisectional curvature if and only if $A$, $B$, $C$ are positive functions of $s$ on $[0,\infty)$.
\end{theorem}

\begin{definition}
	We denote by $\mathcal{M}_n$ the set of all complete $U(n)$ invariant K\"ahler metrics on $\mathbb{C}^n$ with positive bisectional curvature.
\end{definition}

In \cite{K77}, Klembeck constructed an explicit example of a metric in $\mathcal{M}_n$. In \cite{C94,C97}, Cao came up with two examples of K\"ahler Ricci soliton metrics in  $\mathcal{M}_n$. In their paper \cite{WZ11} Wu and Zheng characterized $\mathcal{M}_n$ via a function $\xi= \xi(s)$ and illustrated that the set $\mathcal{M}_n$ is actually quite large.

\begin{definition}[Wu-Zheng \cite{WZ11}] The smooth function $\xi:[0,\infty)\rightarrow \mathbb{R}$ is defined by
	\begin{equation}\label{xi1}
	\xi(s) = -s(\log g'(s))'.
	\end{equation}
\end{definition}

Thus, we have $g'(s) = g'(0) \exp\left(-\displaystyle\int_{0}^{s} \frac{\xi(s)}{s}ds\right)$.

\begin{theorem}[Characterization of $\mathcal{M}_n$ by the function $\xi$, Wu-Zheng \cite{WZ11}]\label{Thm WZ bisec} The metric given by \eqref{M1} is a complete K\"ahler metric with positive bisectional curvature on $\mathbb{C}^n$ if and only if $\xi$ defined by \eqref{xi1} satisfies 
	\begin{equation}\label{xi2}
	\xi(0)= 0, \quad \xi'>0, \quad \xi<1.
	\end{equation}
\end{theorem}
If we let $\Xi$ be the space of all $\xi\in C^\infty[0,\infty)$ satisfying \eqref{xi2}, then $\Xi$ is the space of all diffeomorphisms $[0,\infty)\rightarrow[0,b)$, ($0<b\leq 1$). The space $\Xi$ is in one-to-one corresponence with $\mathcal{M}_n\slash \mathbb{R}^+$.

We will see later that no metric in $\mathcal{M}_n$ satisfies the \textit{extremal condition}.

\subsection{Extremal Condition}\label{section2.3}

For $U(n)$ invariant extremal K\"ahler metrics on $\mathbb{C}^n\backslash \{0\}$, Calabi reduced extremal equation to an ordinary differential equation for the K\"ahler potential $u(s)$, $s=|z|^2$. He then used the function $g(s)=su'(s)$ to write the equation in the form $sg'(s)=F(g(s))$. 

\begin{definition}\label{def extremal}
	We say that a K\"ahler metric satisfies the extremal condition if its scalar curvature $R$ satisfies the Euler equation $R_{,\overline{\alpha} \overline{\beta}}=0$.
\end{definition}

For the rotation invariant K\"ahler metrics on $\mathbb{C}^n \backslash\{0\}$, Calabi \cite{Calabi82} reduced the equation $R_{,\overline{\alpha} \overline{\beta}}=0$ to a nonlinear ordinary differential equation $sg'(s)=F(g(s))$ as follows.

Let us denote by $\G=(\g_{j\overline k})$ the matrix of the K\"ahler metric. Then, as given in \cite{HeLi}, we have 
\begin{align}\label{detG}
\det \G &= (u'(s))^{n-1} (u'(s)+su''(s))\\\label{gjkbar}
\g^{j\overline k} &= e^{v} (u'(s))^{n-2} [(u'(s)+su''(s))\delta_{jk} - u''(s)\overline z_k z_j]
\end{align}
where $v=-\log\det \G$.

Moreover, by direct computation we have
\begin{equation}\label{deldelbarv}
\partial\overline\partial v = -\sum_{j,k=1}^{n} (\delta_{jk} v'(s) +v''(s) \overline z_j z_k) dz_j \wedge d\overline z_k.
\end{equation}
We combine \eqref{gjkbar} and \eqref{deldelbarv} to obtain
\begin{align*}
R &=\sum_{j,k=1}^n \g^{j\overline k} \frac{\partial^2 v}{\partial z_j \partial\overline z_k}\\
& = s^{1-n} e^v [s^n (u'(s))^{n-1}v'(s)]'.
\end{align*}
We substitute the expression for $\det \G$ given in \eqref{detG} into this equation to get
\begin{align}\label{R(s)}
R(s) &= \frac{s^{1-n} [s^n (u')^{n-1}v']'}{(u')^{n-1}(u'+su'')}\\\nonumber
&=\frac{nv'+s(n-1)(u')^{-1} u'' v' +sv''}{u'+su''}\\\nonumber
&=\frac{v'\big(\frac{(n-1)(u'+su'')}{u'} + \frac{u'}{u'}\big) +sv''}{u'+su''}\\\nonumber
&=(n-1)\frac{v'}{u'} + \frac{sv'+s^2v''}{su'+s^2 u''}.
\end{align}
We note that if we substitute $s=e^t$ in \eqref{R(s)}, we obtain Equation (3.9) in \cite{Calabi82}.

The condition that the components of the tangent vector fields $\g^{\alpha\overline\beta} \frac{\partial R}{\partial \overline z_\beta} \partial z_\alpha$ be holomorphic is equivalent to the Euler equation $R_{,\overline\alpha\overline\beta}=0$ (see \cite{Calabi82}). This equation can be expressed in the rotationally symmetric case as follows \cite{Calabi82}:
\begin{equation*}
\g^{\alpha\overline\beta} \frac{\partial R}{\partial \overline z_\beta} = \g^{\alpha\overline\beta} R'(s) z_\beta = z_\alpha\frac{R'(s)}{u'+su''}
\end{equation*}
where, in the last equality, we have used \eqref{detG} and \eqref{gjkbar}. The Euler equation is now equivalent to $\frac{\partial}{\partial\overline z_\beta} \Big(z_\alpha \frac{R'}{u'+su''}\Big) = 0$, \ $\beta=1,\dots,n$; and since $\frac{R'}{u'+su''}$ is real valued, we obtain the equation
\begin{equation*}
\frac{R'}{u'+su''} = \textit{constant}.
\end{equation*}
For convenience, we will set this constant to be $-(n+2)(n+1)c_4$, as in \cite{Calabi82}. We can make use of the variable change $s=e^t$ to integrate the differential equation and obtain
\begin{equation}\label{R}
R=-(n+2)(n+1)c_4 g(s) -(n+1)nc_3
\end{equation}
Replacing $R$ in \eqref{R} by its expression in term of $u$, $v$, and their derivatives, and integrating once more, we obtain Equation (3.12) in Calabi's article \cite{Calabi82}:
\begin{equation}\label{ordinary differential equation}
\frac{g^{n-1} g' }{c_4 g^{n+2} +c_3 g^{n+1}+g^n+c_1g+c_0} = \frac{1}{s}.
\end{equation}

The Euler equation $R_{,\overline\alpha\overline\beta}=0$  has been reduced to an ordinary differential equation $sg'(s)=F(g(s))$ where 
\begin{equation}\label{F(g)}
F(g)= \frac{c_4 g^{n+2} +c_3 g^{n+1}+g^n+c_1g+c_0}{g^{n-1} } .
\end{equation}
After simplification of rational expression in \eqref{F(g)} (if necessary) we will denote the polynomial in the numerator by $H(g)$. If we write $\displaystyle\lim_{s\rightarrow 0^+} g(s)=A$ and $\displaystyle\lim_{s\rightarrow +\infty} g(s)=B\leq \infty$, then we see from Lemma \ref{lemma6.2}, that $H(A)=0$ and $H>0$ on $(A,B)$. Moreover, $H(B)=0$ whenever $B<\infty$.

\subsection{$k$-twisted (Projective) Line Bundles and Orbifolds}\label{k twisted}

Calabi, in his paper \cite{Calabi82}, described $k$-twisted projective line bundles $\mathbb{CP}^1\hookrightarrow \mathcal{F}_k^n \stackrel{\pi}{\twoheadrightarrow} \mathbb{CP}^{n-1} $ for any $k=1,2,\dots$, $n\geq 2$, as follows.

We cover $\mathbb{CP}^{n-1}$ by $n$ coordinate domains $U_\lambda = \{[z_1:\dots:z_n] : z_\lambda \neq 0\}$ ($1\leq \lambda\leq n$). On each $U_\lambda$, we have a holomorphic coordinate system $(_\lambda z^\alpha) = \left(\frac{z_\alpha}{z_\lambda}\right)$, ($1\leq\alpha \leq n$, $\alpha\neq \lambda$). We introduce a projective holomorphic fiber coordinate $y_\lambda \in \mathbb{C}\cup \{\infty\}$ and trivialization $\pi^{-1}(U_\lambda)\simeq U_\lambda\times\mathbb{CP}^1$ ($1\leq\lambda\leq n$) on $\mathcal{F}_k^n$. Here, the transition relation on the fiber coordinate in $\pi^{-1}(U_\lambda\cap U_\mu)$ is given by $$([z_1:\dots:z_n],y_\mu) = \left([z_1:\dots:z_n],\left(\frac{z_\mu}{z_\lambda}\right)y_\lambda\right) $$

We have two distinguished sections $s_0,s_\infty:\mathbb{CP}^{n-1} \rightarrow \mathcal{F}_k^n$ with images denoted by $S_0$ and $S_\infty$, respectively. Here $s_0$ is the zero section given by $y_\lambda=0$, and $s_\infty$ is the infinity section given by $y_\lambda=\infty$ ($1\leq\lambda\leq n$).

We note that the complement $\mathcal{F}_k^n \backslash S_\infty$ gives us the line bundle $\mathcal{O}_{\mathbb{CP}^{n-1}}(-k)$, whereas $\mathcal{F}_k^n \backslash S_0$ gives $\mathcal{O}_{\mathbb{CP}^{n-1}}(k)$. 
Throughout the thesis, we will denote the zero sections of the line bundles $\mathcal{O}_{\mathbb{CP}^{n-1}}(-k)$ and $\mathcal{O}_{\mathbb{CP}^{n-1}}(k)$ by $S_0$ and $S_\infty$, respectively. We will write $\mathring{\mathcal{F}}_k^n$ for the complement $\mathcal{F}_k^n \backslash\{S_0\cup S_\infty\}$. We have a $k:1$ map
\begin{equation}\label{K1}
p:\mathbb{C}^n \backslash\{0\} \rightarrow \mathring{\mathcal{F}}_k^n
\end{equation}
which assigns to any point $(z_1,\dots,z_n)$ with $z_\lambda\neq 0$, the point in $\mathring{\mathcal{F}}_k^n \cap \pi^{-1}(U_\lambda)$ with coordinates $\left( \left(\frac{z_\alpha}{z_\lambda}\right);(z_\lambda)^k\right)$, ($1\leq\alpha \leq n$, $\alpha\neq\lambda$).

The map $p$ induces a biholomorphism 
\begin{equation}\label{K2}
\tilde{p}:\mathbb{C}^n \backslash\{0\} \slash \mathbb{Z}_k \rightarrow  \mathring{\mathcal{F}}_k^n.
\end{equation}
Thus, $\mathcal{O}_{\mathbb{CP}^{n-1}}(-k)$ is obtained by gluing a $\mathbb{CP}^{n-1}$ smoothly into $(\mathbb{C}^n\backslash\{0\})\slash \mathbb{Z}_k $ at $z=0$. Similarly, we obtain  $\mathcal{O}_{\mathbb{CP}^{n-1}}(k)$ if we glue a  $\mathbb{CP}^{n-1}$ smoothly into $(\mathbb{C}^n\backslash\{0\})\slash \mathbb{Z}_k $ at $z=\infty$.

The map $\tilde{p}:\mathbb{C}^n \backslash\{0\} \slash \mathbb{Z}_k \rightarrow   \mathcal{O}_{\mathbb{CP}^{n-1}}(-k) \backslash S_0$ can be written as
$$\tilde{p}(z_1,\dots,z_n) = ([z_1:\dots:z_n];(z_0,\dots,z_n)^{\otimes k}),$$
where $(z_0,\dots,z_n)^{\otimes k}$ denotes the generator of the fiber of $\mathcal{O}_{\mathbb{CP}^{n-1}}(-k)\rightarrow\mathbb{CP}^{n-1}  $ over the point $[z_1:\dots:z_n]$. (See Apostolov-Rollin, \cite{AR17} for more details).

We will denote by $G_k^n$ the compact space obtained from $\mathcal{F}_k^n$ by blowing down its zero section $S_0$ to a point. 
	$G_k^n$ is the weighted projective space $\mathbb{CP}^n_{[k,1,\dots,1]} = (\mathbb{C}^{n+1}\backslash\{0\})/ \sim $ where $(z_0,z_1,\dots,z_n)\sim (\lambda^k z_0, \lambda z_1, \dots, \lambda z_n)$ for any $\lambda \in \mathbb{C}\backslash\{0\}$. $G_k^n$ has an orbifold singularity at $[1:0:\dots:0]$ mordinary differential equationled on $\mathbb{C}^n/\mathbb{Z}_k$.
When $k=1$, $G_k$ is simply $\mathbb{CP}^n$.

\subsection{Closing Conditions}\label{section2.5}
Suppose that we have a $U(n)$-invariant K\"ahler metric $\g$ on $\mathbb{C}^n\backslash\{0\}$, $n\geq 2$, that satisfies the extremal condition $sg'(s)=F(g(s))$, where $g(s)=su'(s)$ $u(s)$ is a K\"ahler potential, and $F(g)$ is given by \eqref{F(g)}.

If  we have a $U(n)$-invariant K\"ahler metric $\g$ on $\mathbb{C}^n\backslash\{0\}$ given  by \eqref{M2}, it induces a metric on $\mathring{\mathcal{F}}_k^n$ via the map $\tilde{p}$. We will denote the induced metric on $\mathring{\mathcal{F}}_k^n$ by $\g$ as well. 

In \cite{Calabi82}, Calabi imposed certain asymptotic conditions on K\"ahler potential $u(s)$ as $s\rightarrow 0^+$ and $s\rightarrow\infty$. These conditions are necessary and sufficient for the metric $\g$  on $\mathring{\mathcal{F}}_k^n$ to be extandable by continuity to a smooth metric on all of ${\mathcal{F}}_k^n$.

Cao \cite{C94} used the map $\tilde p$ to produce $U(n)$-invariant, complete gradient K\"ahler-Ricci soliton (\textit{GKRS}) metrics on line bundles over $\mathbb{CP}^{n-1}$. Feldman-Ilmanen-Knopf \cite{FIK03} generalized this approach by producing $U(n)$-invariant GKRS metrics on $(\mathbb{C}^n\backslash\{0\})\slash \mathbb{Z}_k$, where they allowed new boundary behavior at $z=0$ and $|z|=\infty$. These behaviors are listed as follows.

\begin{enumerate}[label=\arabic*.]
	\item \label{smooth pt}Metric is completed at $z=0$ by adding a smooth point.
	
	\item \label{orb pt} Metric is completed at $z=0$ by adding an orbifold point.
	
	\item \label{smooth cpn at 0}Metric is completed at $z=0$ by adding a smooth or singular $\mathbb{CP}^{n-1}$.
	
	\item \label{complete 0} Metric is complete as $|z|\rightarrow 0$.
	
\end{enumerate}

\begin{enumerate}[label=\alph*.]
	\item \label{smooth cpn} Metric is completed at $z=\infty$ by adding a smooth or singular $\mathbb{CP}^{n-1}$.
	\item \label{complete} Metric is complete as $|z|\rightarrow\infty$.
\end{enumerate}

We note that \ref{smooth pt}\ref{smooth cpn} gives $\mathbb{CP}^n$ (with $k=1$), and \ref{smooth pt}\ref{complete} gives $\mathbb{C}^n$.

The boundary conditions \ref{orb pt}\ref{smooth cpn} correspond to $G_k^n = \mathbb{CP}^n_{[k:1:\dots:1]}$. Conditions \ref{smooth cpn at 0}\ref{complete} give  $\mathring{\mathcal{F}}_k^n \cup S_0=  \mathcal{O}_{\mathbb{CP}^{n-1}}(-k)$, and conditions \ref{complete 0}\ref{smooth cpn} give  $\mathring{\mathcal{F}}_k^n \cup S_\infty=  \mathcal{O}_{\mathbb{CP}^{n-1}}(k)$. Calabi's compact $k$-twisted $\mathbb{CP}^{1}$-bundle ${\mathcal{F}}_k^n $ is obtained via \ref{smooth cpn at 0}\ref{smooth cpn}.

A $U(n)$ invariant K\"ahler metric on $\mathring{\mathcal{F}}_k^n$ induced by the map \eqref{K2} cannot be completed by adding a $\mathbb{CP}^{n-1}$ at $z=0$ and a smooth or orbifold point at $z=\infty$. This follows since $g(s)$ is a strictly increasing function of $s$.

Adding a $\mathbb{CP}^{n-1}$ smoothly to $(\mathbb{C}^n\backslash \{0\}) \slash \mathbb{Z}_k$ corresponds to a simple zero of $F$ \cite{Calabi82}.  In what follows,  we explain this correspondence as it is presented in \cite{FIK03}. 

If the sign of $F'$ at the simple root is positive (resp. negative), it means we are adding $\mathbb{CP}^{n-1}$ at $z=0$ (resp. $|z|=\infty$). This can be seen as follows. Let $\displaystyle\lim_{s\rightarrow 0^+}g(s)=A$, $\displaystyle\lim_{s\rightarrow \infty}g(s)=B$ ($0<A<B\leq \infty$). By Lemma \ref{lemma6.1}, we have $H(A)=0$ and $H>0$ on $(A,B)$. If $A>0$ is a simple root of $F$, then it is a simple root of $H$. In this case, we have $H'(A)>0$, implying $F'(A)>0$. Similarly, if $B<\infty$ is a simple root, we have $F'(B)<0$. Therefore, the sign of $F'$ at a simple root determines whether $\mathbb{CP}^{n-1}$ is added at $z=0$ or at $|z|=\infty$.

Let us assume
\begin{equation}\label{C1}
F(A)=0,\quad A>0,\quad F'(A)=\theta>0.
\end{equation}
For convenience, we will switch to a new parameter $t=\log s$, $-\infty<t<\infty$, as in \cite{FIK03}. We will obtain a specific form for $g(s)$ in a neighborhood of $s=0$.

We write the ordinary differential equation $sg'(s)=F(g(s))$ in the form $\phi'(t)=F(\phi(t))$, where $\phi(t):=g(s)$. We have $\phi'(t)=sg'(s)>0$, hence $t=t(\phi)$ is a smooth strictly increasing function of $\phi$. We have a diffeomorphism $\psi=\Phi(\phi)$, given by $\Phi(\phi):=e^{\theta t(\phi)}$. The ordinary differential equation $\phi'(t)=F(\phi(t))$ is conjugate to the equation $\psi'(t)=\theta \psi(t)$, so $\phi(t)$ has the form 
$$\phi(t)=A + e^{\theta t } G_0(e^{\theta t})$$ 
as $t\rightarrow -\infty$. Here $G_0$ is a smooth function on $(-\epsilon,\epsilon)$ with $G_0(0)>0$ \cite{Calabi82,C94,FIK03}.

Let us switch back to the parameter $s=e^t$. We have seen that \eqref{C1} implies 
\begin{equation}\label{C2}
g(s)=A+s^{\theta} G_0(s^\theta)
\end{equation}
where $G_0$ is given as above. Since $F'(A)=\theta>0$, we are adding $\mathbb{CP}^{n-1}$ at $z=0$. It follows from \eqref{M2} and \eqref{K2} that each complex line through the origin has a cone angle $2\pi\theta/k$.

\begin{remark}\label{remark geod dist}
	When we have $\displaystyle\lim_{s\rightarrow 0^+}g(s)=A>0$, $F(A)=0$, $F'(A)>0$, equation \eqref{C2} implies that geodesic distance to $z=0$ is finite, i.e. $\displaystyle\frac{1}{2} \int_{0}^{s_0} \sqrt{\frac{g'(s)}{s}}\ ds <\infty$.
\end{remark}

A similar discussion follows when we have 
\begin{equation}\label{C3}
F(B)=0, \quad B>0,\quad F'(B)=-\theta<0.
\end{equation}
Let us assume $g(s)$ solves $sg'(s)=F(g(s))$, with $\displaystyle\lim_{s\rightarrow \infty}g(s)=B$ and \eqref{C3} is satisfied. Then we have
\begin{equation}\label{C4}
g(s) = B+ s^{-\theta} G_\infty (s^{-\theta})
\end{equation}
where $G_\infty$ is smooth on $(-\epsilon,\epsilon)$ and $G_\infty(0)<0$.

For the proof of the following Lemma, see \cite{FIK03}.
\begin{lemma}[Calabi \cite{Calabi82}]\label{Lemma Calabi} Let $n\geq 2$.
	\begin{enumerate}
		\item When $\theta=k$ in \eqref{C2}, the induced K\"ahler metric is smooth on a neighborhood of the zero section in $\mathcal{O}_{\mathbb{CP}^{n-1}} (-k)$.
		
		\item When $\theta=-k$ in \eqref{C4}, the induced K\"ahler metric is smooth on a neighborhood of the zero section in $\mathcal{O}_{\mathbb{CP}^{n-1}} (k)$.
	\end{enumerate}
\end{lemma}

\begin{remark}
	In Section \ref{section un inv on cn} Lemma \ref{addingsmoothpt}, we will see that a solution of $sg'(s)=F(g(s))$ on $\mathbb{C}^n\backslash\{0\}$ gives a smooth metric on $\mathbb{C}^n$ if and only if $F(0)=0$. In this case, we say that we are adding a smooth point at $z=0$.
\end{remark}

\section{$U(n)$ invariant Extremal K\"ahler Metrics on $\mathbb{C}^n$}\label{section un inv on cn}

\subsection{List of Solutions on $\mathbb{C}^n$ and Related Results}\label{section 3.1}

We start with a smooth extension lemma, which, together with Theorem \ref{n smooth solution} and Theorem \ref{Thm2.2} will give rise to an improved version (Lemma \ref{addingsmoothpt}).

\begin{lemma}\label{starlemma}
	Let $u\in C^\infty (0,\infty)$ be the potential of a rotation invariant K\"ahler metric on $\mathbb{C}^n\backslash\{0\}$ that satisfies the extremal condition. Then, the metric extends smoothly to $\mathbb{C}^n$ if and only if $u\in C^2 [0,\infty)$ and $\lim\limits_{s\rightarrow 0^+} u'(s)$ is positive.
\end{lemma}
\begin{proof}
	If we assume $u\in C^2[0,\infty)$, then we have $\displaystyle\lim_{s\rightarrow 0^+} g(s) = \lim_{s\rightarrow 0^+} su'(s)=0$, and $\displaystyle\lim_{s\rightarrow 0^+} sg'(s)=0$. It follows from Equation \eqref{ordinary differential equation} that $c_0=c_1=0$. Equation \eqref{ordinary differential equation} can be written as 
	\begin{align}\label{star}
	u'' = c_4 s(u')^3 + c_3 (u')^2.
	\end{align}
	Differentiating \eqref{star}, we obtain $u\in C^{\infty} [0,\infty)$. This implies that the hypothesis of Monn's smoothness result\footnote{One has to make a parameter change $r=\sqrt s$, and use Proposition \ref{MonnProp4.1}.} (see Proposition \ref{MonnProp2.1} below) for the corresponding radial function $u(z_1,\dots,z_n)$ is satisfied for all $k\geq 0$. Together with the condition $\lim\limits_{s\rightarrow 0^+} u'(s)>0$, we  conclude that the metric extends smoothly to $\mathbb{C}^n$.
	
	The converse is clear.
\end{proof}

We state He-Li's Theorem 1.1 in \cite{HeLi} that gives a complete list of rotation invariant cscK metrics on $\mathbb{C}^n$.

\begin{theorem}[\cite{HeLi}, Theorem 1.1]\label{n smooth solution}
	Suppose $n \geq 2$ is an integer.
	\begin{enumerate}
		\item \label{metric std} The rotation invariant K\"ahler metric $\omega$ with zero constant scalar curvature on $\mathbb{C}^n$ must be a multiple of the standard Euclidean metric.
		
		\item The rotation invariant K\"ahler metric $\omega$ with constant scalar curvature $n(n+1)$ on $\mathbb{C}^n$ must be of the form
		\[
		\omega=i\frac{\sum_{j=1}^{n}dz_j \wedge d \overline{z}_j}{\sum_{j=1}^{n}|z_j|^2+a}-i\frac{(\sum_{j=1}^{n}\overline{z}_jdz_j)\wedge \overline{(\sum_{j=1}^{n}\overline{z}_jdz_j)}}{(\sum_{j=1}^{n}|z_j|^2+a)^2}
		\]
		where $a>0$ is a constant.
		
		\item There does not exist rotation invariant K\"ahler metric with negative constant scalar curvature on $\mathbb{C}^n$.
	\end{enumerate}
\end{theorem}

The following theorem gives a complete list of solutions of Calabi's $U(n)$-invariant extremal equation on $\mathbb{C}^n$. We note that for all real numbers $\alpha$, $\beta$, $\gamma$, $c$ satisfying given conditions, a unique non-degenerate metric defined for all $s\in[0,\infty)$ does exist.

{
	\renewcommand{\thetheorem}{\ref{Thm2.2}}
	\begin{theorem}
		Let $n\geq 2$ and $u:[0,\infty)\rightarrow \mathbb R$ be 
		the potential of an extremal K\"ahler metric 
		on $\mathbb C^n$. Then, for $g(s) = su'(s)$, one of the following is true:
		\begin{enumerate}
			\item $\omega$ is a cscK metric.
			\item There exist constants $\beta,c$ with $\beta>0$ such that $g(s)$ is the smooth strictly increasing function $g:(0,\infty) \rightarrow(0,\beta)$ determined by the ordinary differential equation
			$$sg'(s)=F(g(s)) = \frac{1}{\beta^2}g(s) (g(s)-\beta)^2.$$ 
			Integrating, we get
			$$\log(g(s)) -\log(\beta-g(s))-\beta\frac{1}{g(s)-\beta} =\log s+c.$$
			These metrics have positive nonconstant scalar curvature.
			
			\item There exist constants $\gamma,\beta,c$ with $\gamma<0<\beta$ such that  $g(s)$ is the smooth strictly increasing function ranging from $g:(0,\infty)\rightarrow(0, \beta)$ determined by the ordinary differential equation
			$$sg'(s)=F(g(s)) = \frac{1}{\beta\gamma}g(s) (g(s)-\beta)(g(s)-\gamma).$$
			Integrating, we get 
			$$\log(g(s))+\frac{\beta\gamma}{\beta(\beta-\gamma)}\log(\beta - g(s)) +\frac{\beta\gamma}{\gamma(\gamma-\beta)}\log(g(s)-\gamma) = \log(s)+c.$$
			
			\item There exist constants $\gamma,\beta,c$ with $0<\beta<\gamma$ such that  $g(s)$ is the smooth strictly increasing function  $g:(0,\infty)\rightarrow(0, \beta)$ determined by the ordinary differential equation
			$$sg'(s)=F(g(s)) = \frac{1}{\beta\gamma}g(s) (g(s)-\beta)(g(s)-\gamma).$$
			Integrating, we get
			$$\log(g(s))+\frac{\beta\gamma}{\beta(\beta-\gamma)}\log(\beta - g(s)) +\frac{\beta\gamma}{\gamma(\gamma-\beta)}\log(\gamma-g(s)) = \log(s)+c.$$
		\end{enumerate}
	\end{theorem}
	\addtocounter{theorem}{-1}
}

\begin{proof} 
	From the proof of Lemma \ref{starlemma}, if $u$ is the potential of a smooth metric on $\mathbb{C}^n$, then we have $\displaystyle\lim_{s\rightarrow 0^+} g(s)=\displaystyle\lim_{s\rightarrow 0^+} su'(s)=0$ and the constants $c_0$ and $c_1$ in equation \eqref{ordinary differential equation} vanish. Equation \eqref{ordinary differential equation} becomes
		\begin{equation}\label{odesmooth}
		sg'(s) = F(g(s)) = c_4 g^{3} + c_3 g^{2} + g 
		\end{equation}
	In this case, the polynomial $H(x)$ in Lemma \ref{lemma6.2} is given by $c_4x^3+c_3x^2+x$, and unless $c_3=c_4=0$, we have $B<\infty$ for degree reasons. By the same lemma, $H(A)=0$ and $H>0$ on $(A,B)$, and all roots are real. We have $A=\displaystyle\lim_{s\rightarrow 0^+} g(s)=0$. 
	
	\begin{enumerate}[label=\underline{\smash{\textit{Case}\,(\arabic*)}}, align=left, leftmargin=0pt]
		\item \textit{$c_4=0$} \\We see from equation \eqref{R} that $\omega$ is a cscK metric.\\
		
		\item \textit{$c_4\neq 0$ and the polynomial $H(x)=c_4x^3+c_3x^2+x$ has roots $\alpha,\beta,\beta$ with $\alpha<\beta$}\\
		It follows from Lemma \ref{lemma6.2} that $\alpha=A=0$. Since $H(x)>0$ in $(\alpha,\beta)$, we have $c_4>0$, and the equation \eqref{ordinary differential equationsmooth} can be written as 
		$$g'\beta^2 \left\{\frac{1}{\beta^2} \frac{1}{g(s)}-\frac{1}{\beta^2}\frac{1}{g(s)-\beta} +\frac{1}{\beta} \frac{1}{(g(s)-\beta)^2} \right\} = \frac{1}{s}.$$
		There exists a constant $c$ such that
		\begin{equation}\label{11}
		\log(g(s))-\log(\beta-g(s)) -\frac{\beta}{g(s)-\beta} = \log s+ c.
		\end{equation} 
		Since $H(x)=c_4x^3+c_3 x^2 +x = c_4 x(x-\beta)^2$,  we have $c_4\beta^2=1$.
		
		On the other hand, Lemma \ref{lemma6.1} implies that there exists a unique smooth strictly increasing function $g(s)=su'$, $g:(0,\infty)\rightarrow(0,\beta)$ satisfying \eqref{11}. We note that this case yields solution \eqref{metric beta} in the Theorem. \\
		
		\item \textit{$c_4\neq 0$ and the polynomial $H(x)=c_4x^3+c_3x^2+x$ has roots $\alpha,\alpha,\beta$ with $\alpha<\beta$}\\
		It follows from Lemma \ref{lemma6.2} that $\alpha=A=0$, but obviously this polynomial does not have a double root at $0$.\\
		
		\item \textit{$c_4\neq 0$ and the polynomial $H(x)=c_4x^3+c_3x^2+x$ has distinct roots $\gamma<\beta<\alpha=0$}\\
		It follows from Lemma \ref{lemma6.2} that $B<\infty$. Therefore $B$ must be a root of $H$, and we must have $H>0$ on $(0,B)$. This contradicts our choice of roots $\gamma<\beta<\alpha =0$.\\
		
		\item \textit{$c_4\neq 0$ and the polynomial $H(x)=c_4x^3+c_3x^2+x$ has distinct roots $\gamma<\alpha=0<\beta$}
		
		By Lemma \ref{lemma6.2} we have $A=\alpha=0$, $B=\beta$, and $H(x)>0$ on $(0,\beta)$. Then $c_4<0$, $c_4\beta\gamma=1$, and the equation \eqref{ordinary differential equationsmooth} can be written as
		\begin{equation}\label{12}
		g'\left\{\frac{1}{g(s)}+\frac{\beta\gamma}{\beta(\beta-\gamma)}\frac{1}{g(s)-\beta} +\frac{\beta\gamma}{\gamma(\gamma-\beta)}\frac{1}{g(s)-\gamma} \right\} = \frac{1}{s}.
		\end{equation}
		There exists a constant $c$ such that 
		$$\log(g(s))+\frac{\beta\gamma}{\beta(\beta-\gamma)}\log(\beta-g(s)) + \frac{\beta\gamma}{\gamma(\gamma-\beta)}\log(g(s)-\gamma)=\log s +c. $$
		It follows from Lemma \ref{lemma6.1} that there exists a unique smooth strictly increasing function $g(s):(0,\infty)\rightarrow(0,\beta)$. We note that this case yields solution \eqref{metric betagamma} in the Theorem.\\
		
		\item \textit{$c_4\neq 0$ and the polynomial $H(x)=c_4x^3+c_3x^2+x$ has distinct roots $\alpha=0<\beta<\gamma$}\\
		It follows from Lemma \ref{lemma6.2} that $A=\alpha=0$, $B=\beta$ and $c_4>0$. As in the previous case, Equation \eqref{ordinary differential equationsmooth} can be rewritten as Equation \eqref{12}. Integrating both sides of \eqref{12} we get
		$$\log(g(s))+\frac{\beta\gamma}{\beta(\beta-\gamma)}\log(\beta-g(s)) + \frac{\beta\gamma}{\gamma(\gamma-\beta)}\log(\gamma-g(s))=\log s +c. $$
		If we let $h(t):(0,\beta)\rightarrow \mathbb{R}$ be the function 
		$$\log(t)+\frac{\beta\gamma}{\beta(\beta-\gamma)}\log(\beta-t) + \frac{\beta\gamma}{\gamma(\gamma-\beta)}\log(\gamma-t) $$
		we see that $\displaystyle \lim_{t\rightarrow 0^+} h(t)=-\infty$,  $\displaystyle \lim_{t\rightarrow \beta^-} h(t)=\infty$, and $h'(t)>0$ on $(0,\beta)$.  Lemma \ref{lemma6.1} guarantees the unique existence of a smooth $g(s):(0,\infty)\rightarrow(0,\beta)$ with the desired properties. We note that this case yields solution \eqref{metric betagamma2} in the Theorem.
	\end{enumerate}
\end{proof}

\begin{corollary}\label{thm1}
	There does not exist a rotation invariant extremal K\"ahler metric with negative scalar curvature   on $\mathbb{C}^n$.
\end{corollary}
\begin{proof}
	See Section \ref{section proof}.
\end{proof}

\begin{remark}\label{remark extension}
	We note that, in the proof of Theorem \ref{Thm2.2}, in order to obtain the implicit solutions given by cases (1)--(4), we have not used the full strength of $u(s)$ being in $C^{\infty}[0,\infty)$. In the proof, we have only used $u(s)\in C^{\infty}(0,\infty)$, $\lim\limits_{s\rightarrow 0^+} g(s)=0$, and $c_0=c_1=0$. A careful inspection of the implicit solutions (1)--(4) shows that $\lim\limits_{s\rightarrow 0^+} u'(s)$ is finite and positive, hence, equation \eqref{star} implies that $u(s)$ is in $C^{\infty}[0,\infty)$. It follows from Lemma \ref{starlemma} that such metrics can be smoothly extended to the origin.
\end{remark}

The following lemma tells us when a rotation invariant K\"ahler metric with extremal condition on $\mathbb{C}^n\backslash\{0\}$ can be smoothly extended to $\mathbb{C}^n$, $n\geq 2$.

\begin{lemma}[Adding a smooth point at $z=0$]\label{addingsmoothpt} Let $g:(0,\infty)\rightarrow (A,B)$ ($0\leq A<B\leq \infty$) be a positive, strictly increasing solution of $sg'=F(g(s))$. Then the following are equivalent.
	\begin{enumerate}
		\item $g$ induces a smooth extremal K\"ahler metric on $\mathbb{C}^n$.
		\item $\displaystyle \lim_{s\rightarrow 0^+} g(s)=0$.
		\item $F(0)=0$.
		\item $g$ induces a Bochner-K\"ahler metric on $\mathbb{C}^n$.
	\end{enumerate}	
\end{lemma}

\begin{proof}
	See Section \ref{section proof}.
\end{proof}

Now that we have a complete list of $U(n)$ invariant extremal K\"ahler metrics on $\mathbb{C}^n$, we can utilize the formulation given in \cite{WZ11} for complete $U(n)$ invariant K\"ahler metrics with positive bisectional curvature.

We recall that  the set of all complete $U(n)$ invariant K\"ahler metrics on $\mathbb{C}^n$ with positive bisectional curvature is denoted by $\mathcal{M}_n$, as given in \cite{WZ11}.

\begin{proposition} \label{prop2.8}
	There is no complete $U(n)$ invariant extremal K\"ahler metric on $\mathbb{C}^n$ with positive bisectional curvature.
\end{proposition}
\begin{proof}
	We have only two types of complete $U(n)$ invariant extremal K\"ahler metrics on $\mathbb{C}^n$. The first type is given by \eqref{metric std} of Theorem \ref{n smooth solution}, namely a scalar multiple of the standard Euclidean metric on $\mathbb{C}^n$. Metrics of this type clearly do not have positive bisectional curvature.
	
	The second type is given by \eqref{metric beta} of Theorem \ref{Thm2.2}. See Example \ref{biquard} below for a proof that this metric is complete. We note that all cases other than these two in Theorem \ref{n smooth solution} and Theorem \ref{Thm2.2} give incomplete metrics. This follows from Remark \ref{remark geod dist}; if $\displaystyle\lim_{s\rightarrow +\infty} F(g(s))<0$  is finite, then the corresponding metric is incomplete on $\mathbb{C}^n$. We will see that bisectional curvature is not positive in this case either. We will compute the $\xi$ function for this metric, and show that it does not satisfy the  properties given in Theorem \ref{Thm WZ bisec}.
	
	By definition we have $\xi=-s(\log(g'(s)))'$. We recall that the ordinary differential equation $sg'(s) = F(g(s))$ is given by 
	$$sg'(s) = \frac{1}{\beta^2}g(s)(g(s)-\beta)^2.$$
	We compute 
	\begin{align*}
	(\log (g'))' &=-\frac{1}{s} + \frac{(g(s)-\beta)^2}{\beta^2 s} - \frac{2g(s)(\beta-g(s))}{\beta^2 s}.
	\end{align*}
	Then we have 
	\begin{align*}
	\xi(s)&=1-\frac{1}{\beta^2}(g(s)-\beta)^2 + \frac{2}{\beta^2}g(s)(\beta-g(s))\\
	&= -\frac{1}{\beta^2}g (3g-4\beta)
	\end{align*}
	$\xi(s)$ is a polynomial in $g$ restricted to the interval $(0,\beta)\ni g$. We see that $\displaystyle\frac{d\xi}{ds} = \frac{d\xi}{dg}\frac{dg}{ds}$ is not positive on $(0,\infty)\ni s$. Therefore $\xi$ fails to satisfy the necessary and sufficient conditions in Theorem \ref{Thm WZ bisec}. The metric in Case \eqref{metric beta} of Theorem \ref{Thm2.2} does not have positive bisectional curvature.
	
\end{proof}

\begin{remark}\label{rmrk-Bryant} 
	In \cite{B01}, Bryant shows that complete Bochner-K\"ahler metrics on $\mathbb C^n$ with non-constant scalar curvature given by Tachibana-Liu ansatz are unique up to constant multiples and scaling. Bryant's result is reflected in our discussion in the proof of Proposition \ref{prop2.8}, where we show that there is no complete $U(n)$ invariant extremal K\"ahler metric on $\mathbb C^n$  with positive bisectional curvature.
\end{remark}

\subsection{Examples of Extremal K\"ahler Metrics on $\mathbb{C}^n$}\label{section examples}

He-Li's Theorem \ref{n smooth solution} and our Theorem \ref{Thm2.2} together give a complete list of solutions for Calabi's $U(n)$ invariant extremal equation on $\mathbb{C}^n$. In this section, we retrieve the induced complete extremal K\"ahler metrics on $\mathbb{C}^n$, $\mathbb{CP}^n$, and $\mathbb{CP}_{[k,1,\dots,1]}^n$. We note that these metrics are Bochner-K\"ahler.

\begin{remark}\label{rmk-BKL}
	The solutions given by Theorem \ref{Thm2.2} induce the following metrics.
	\begin{enumerate}[\itshape(i)]
		\item Case \ref{cscK} is the constant scalar curvature case. It induces multiples of Euclidean and Fubini-Study metrics.

		\item Case \ref{metric beta} induces a complete extremal metric on $\mathbb{C}^n$ with non-constant scalar curvature.  (See Example \ref{biquard}). This metric is due to Olivier Biquard.
		
		\item Case \ref{metric betagamma} induces strictly extremal metrics on $\mathbb{CP}^n_{[k,1,\dots,1]}$. (See Example \ref{ex gkn}). These metrics were studied in various settings. They were referred to as degenerate extremal K\"ahler metrics on Hirzebruch surfaces in \cite{HS97}, as self-dual K\"ahler metrics in \cite{AG02} ($n=2$), and as Bochner-K\"ahler metrics in \cite{B01}. See \cite{Abreu01} for a symplectic approach. 
		
		\item Case \ref{metric betagamma2} induces strictly extremal metrics on $(\mathbb{CP}^n, \mathbb{CP}^{n-1})$. (See Example \ref{ex gaud}). These metrics first appeared in \cite{G19}.
	\end{enumerate}
\end{remark}

Case  \eqref{metric beta} of Theorem \ref{Thm2.2} gives a complete extremal metric on $\mathbb{C}^n$ with non-constant scalar curvature. This metric was first discovered by Olivier Biquard. In the next example, we obtain this metric.

\begin{example}[A complete $U(n)$ invariant extremal K\"ahler metric on $\mathbb{C}^n$] \label{biquard}
	We will see that \eqref{metric beta} of Theorem \ref{Thm2.2} induces complete metrics on $\mathbb{C}^n$. In this case, extremal equation is given by $sg'(s)=F(g(s))$ where 
	$$F(g)=c_4 g^3 + c_3 g^2 +g =c_4 g(g-\beta)^2.$$
	Here we have $A=\displaystyle\lim_{s\rightarrow 0^+}g(s)=0$, $B=\displaystyle\lim_{s\rightarrow +\infty}g(s)=\beta<\infty$, and   $c_4 =\frac{1}{\beta^2}>0$.
	
	We will show that geodesic distance from a point $z_0$ to $|z|=\infty$ is infinite, i.e. 
	$$\int_{s_0}^{\infty}\sqrt{\frac{g'(s)}{s}}\ ds = \int_{s_0}^{\infty} \frac{\sqrt{F(g(s))}}{s}\ ds=\infty.$$
	There exists a $d_1>0$ such that on $(s_0,\infty)$ we have 	
	$$\sqrt{F(g(s))}= \frac{1}{\beta}|g-\beta| \sqrt{g} >d_1 (\beta-g)$$
	and 
	$$\int_{s_0}^{\infty} \frac{\sqrt{F(g(s))}}{s}\ ds >d_1 \int_{s_0}^{\infty} \frac{\beta-g}{s}\ ds.$$
	The solution $g(s)$ is given by Equation \eqref{11} as follows. 
	$$\log(g(s))-\log(\beta-g(s)) + \frac{\beta}{\beta-g(s)} = \log s+ c$$
	The term $\log(g(s))$ is bounded on $(s_0,\infty)$. We can choose $s_0$ large enough so that $\log(\beta-g(s))<0$ and $\log s - \log(g(s))+ c>0$ on $(s_0,\infty)$. In this case, Equation \eqref{11} implies $\frac{\beta}{\beta-g(s)}<\log s + c_1$. Therefore
	$$\int_{s_0}^{\infty} \sqrt{\frac{g'(s)}{s}}\ ds >d_1 \int_{s_0}^{\infty} \frac{\beta-g(s)}{s}\ ds > \beta\, d_1 \int_{s_0}^{\infty}\frac{ds}{s(\log s+c_1)} = \infty.$$
	The metric is complete on $\mathbb{C}^n$.
	
	Given that $0<g(s)<\beta$, $c_4 = \frac{1}{\beta}^2$, $c_3=-\frac{2}{\beta}$, and $R=-(n+2)(n+1)c_4g(s) - (n+1)n c_3$, we conclude that $R>0$ for all $s\in[0,\infty)$, hence on $\mathbb{C}^n$.
\end{example}

The positive line bundle $\mathcal O_{\mathbb{CP}^{n-1}}(k)$, $k=1,2,\dots$ is obtained by gluing a $\mathbb{CP}^{n-1}$ to $(\mathbb{C}^{n}\backslash \{0\})\slash \mathbb{Z}_k$ at $|z|=\infty$. If we compactify $\mathcal O_{\mathbb{CP}^{n-1}}(k)$ by adding a singular point at $z=0$, we obtain the weighted complex projective space $G_k^n = \mathbb{CP}^n_{[k,1,\dots,1]}$. Alternatively, we can obtain $G_k^n$ by blowing down $|z|=0$ section in $\mathcal{F}_k^n$. We note that the singular point $z=0$ is modeled on $\mathbb{C}^n\slash \mathbb{Z}_k$. Here, we show that Case \eqref{metric betagamma} of Theorem \ref{Thm2.2} gives an extremal K\"ahler metric with \emph{non-constant scalar curvature} on $G_k^n$ ($n\geq 2$).

\begin{example}[Strictly extremal metrics on $G_k^n$, $n\geq 2$]\label{ex gkn}
		Let us consider Case \eqref{metric betagamma} of Theorem \ref{Thm2.2}. Since $c_4\neq 0$, it follows from Equation \eqref{R} that this has non-constant scalar curvature.
		The ordinary differential equation  is given by $sg'(s)=F(g(s)) = \frac{1}{\beta\gamma} g(s) (g(s)-\beta)(g(s)-\gamma)$. Here, we have $\gamma<0<\beta$, $\displaystyle\lim_{s\rightarrow 0^+}g(s)=0$, $\displaystyle\lim_{s\rightarrow +\infty}g(s)=\beta$.
		
		A $U(n)$ invariant K\"ahler metric on $\mathbb{C}^n $ induces a smooth orbifold metric on $G_k^n\backslash S_\infty$ via the $k:1$ map $p$ given by \eqref{K1}. Here, $S_\infty$ stands for the zero section of $\mathcal O_{\mathbb{CP}^{n-1}}(k)$. It follows from Lemma \ref{Lemma Calabi} that the induced metric can be extended smoothly to $S_\infty$ if and only if $F(\beta)=0$ and $F'(\beta)=-k$. 
		
		We clearly have $F(0)=F(\beta)=0$. We compute $F'(\beta)=\frac{\beta-\gamma}{\gamma}$. For every integer $k\geq 2$, there exist $\gamma,\beta$ ($\gamma<0<\beta$) that satisfy $F'(\beta)=-k$. Namely, let $\beta=|\gamma|(k-1)$. This gives a 1-parameter family of extremal metrics on $G_k^n$, for $k,n\geq 2$.
		
\end{example}

In \cite{DL16}, Dabkowski-Lock explicitly constructed a family of extremal K\"ahler edge cone metrics on $(\mathbb{CP}^2,\mathbb{CP}^{1})$ with cone angles $2\pi\theta$, $\theta\geq 0$. Here, we give examples of strictly extremal metrics on $(\mathbb{CP}^n,\mathbb{CP}^{n-1})$ with cone angles $2\pi\theta$, $0<\theta<1$, $n\geq 2$. These metrics first appeared in \cite{G19}.

\begin{example}  \textbf{(Strictly extremal metrics on  $(\mathbb{CP}^n,\mathbb{CP}^{n-1})$ with cone angles $2\pi\theta$, $0<\theta<1$)} \label{ex gaud}
		Let us consider Case \eqref{metric betagamma2} of Theorem \ref{Thm2.2}. Since $c_4\neq 0$, it follows from Equation \eqref{R} that this is an extremal metric with non-constant scalar curvature on $\mathbb{C}^n$. Extremal equation is given by  $sg'(s)=F(g(s)) = \frac{1}{\beta\gamma} g(s) (g(s)-\beta)(g(s)-\gamma)$. Here, we have $0<\beta<\gamma$, $\displaystyle\lim_{s\rightarrow 0^+}g(s)=0$, $\displaystyle\lim_{s\rightarrow +\infty}g(s)=\beta$. 
		
		As in the previous example we compute $F(0)=F(\beta)=0$ and $F'(\beta)=\frac{\beta-\gamma}{\gamma}=-\theta$. The inequality  $0<\beta<\gamma$ implies that $0<\theta<1$. Therefore, we get a family of metrics determined by $\beta$, $\gamma$ ($0<\beta<\gamma$) and $\theta=1\frac{\beta}{\gamma}$. These are extremal metrics on $(\mathbb{CP}^n,\mathbb{CP}^{n-1})$ with non-constant positive scalar curvature and with cone angle $2\pi\theta$, $0<\theta<1$, along $\mathbb{CP}^{n-1}$ attached at $|z|=\infty$.
\end{example}

\subsection{Proofs}\label{section proof} In this section we give the proofs of Lemma \ref{addingsmoothpt} and Corollary \ref{thm1}.

\begin{proof}[Proof of Lemma \ref{addingsmoothpt}]   
	
	\underline{{(1) $\Rightarrow$ (2)}} Let us assume we have a smooth $U(n)$ invariant metric on $\mathbb{C}^n$. Then we have $u(s)\in C^{\infty}[0,\infty)$. This implies $\displaystyle \lim_{s\rightarrow 0^+} g(s)=\lim_{s\rightarrow 0^+} su'(s)= 0$.
	
	\underline{{(2) $\Rightarrow$ (1)}} We will see that $\displaystyle \lim_{s\rightarrow 0^+} g(s)=0 $ is a sufficient condition for the metric to extend smoothly to $z=0$. $\displaystyle \lim_{s\rightarrow 0^+} g(s)=0$ implies that the constants $c_0$ and $c_1$ in equation \eqref{ordinary differential equation} vanish. This can be seen as follows. Assume  $\displaystyle \lim_{s\rightarrow 0^+} g(s)=A=0$ and $c_0\neq 0$. The condition $c_0\neq 0$ implies
	$$F(g)=\frac{c_4 g^{n+2} + c_3 g^{n+1} + g^n + c_1 g  +c_0}{g^{n-1}} = \frac{H(g)}{g^{n-1}}$$
	and $H(0)=c_0\neq 0$.  This contradicts (1) of Lemma \ref{lemma6.2} which requires $H(A)=0$. So we must have $c_0=0$. 
	
	Now let us assume $\displaystyle \lim_{s\rightarrow 0^+} g(s)=0$, $c_0=0$, and $c_1\neq 0$. Then 
	$$F(g)=\frac{c_4 g^{n+1} + c_3 g^{n} + g^{n-1} + c_1 }{g^{n-2}} = \frac{H(g)}{g^{n-2}}$$
	and $H(0)=c_1\neq 0$, which contradicts (1) of Lemma \ref{lemma6.2} again. Therefore $\displaystyle \lim_{s\rightarrow 0^+} g(s)=0$ implies $c_0=c_1=0$.
	
	It follows from Remark \ref{remark extension} that, since we have $\displaystyle \lim_{s\rightarrow 0^+} g(s)=0$ and $c_0=c_1=0$,  the metric smoothly extends to the origin.
	
	\underline{{(2) $\Rightarrow$ (3)}} Let us assume $\displaystyle \lim_{s\rightarrow 0^+} g(s)=0$. We have already seen that this implies $c_0=c_1=0$. It follows from the definition of $F$ that $F(0)=0$.
	
	\underline{{(3) $\Rightarrow$ (2)}} $F(0)=0$ implies $c_0=c_1=0$. This can be seen from the definition of $F$ ($n\geq 2$) and the limit 
	$$\lim\limits_{x\rightarrow 0} \frac{c_1 x+c_0}{x^{n-1}} =\lim_{x\rightarrow 0} (F(x) -c_4 x^3 -c_3 x^2 -x) = 0.$$
	Now, we will show that $c_0=c_1=0$ implies $\displaystyle \lim_{s\rightarrow 0^+} g(s)=0$.
	
	Let us assume $\displaystyle \lim_{s\rightarrow 0^+} g(s)=A>0$. We will arrive at a contradiction. If $c_0=c_1=0$, equation \eqref{ordinary differential equation} becomes $sg'= c_4g^3 + c_3 g^2 + g = H(g)$.
	
	We have the following cases:
	\begin{itemize}
		\item \fbox{$c_4=c_3=0$}\\ In this case $H(g)=g$ and $H(A)\neq 0$ for $A>0$. This contradicts (1) of Lemma \ref{lemma6.2}. 
		
		\item \fbox{$c_4=0$, $c_3\neq 0$}\\ We have $H(g) = g(c_3 g + 1)$. Since $A>0$ and $H(A)$ vanishes by (1) of Lemma \ref{lemma6.2}, we have $B=\infty$. But this contradicts (2) of Lemma \ref{lemma6.2} for degree reasons.
		
		\item \fbox{$c_4\neq0$}\\ It follows from (2) of Lemma \ref{lemma6.2} that $B<\infty$. We have $$H(g)=c_4 g^3 + c_3 g^2 +g = c_4 g(g-A)(g-B)$$ and $H>0$ on $(A,B)$, ($0<A<B<\infty$). This implies $c_4<0$, which contradicts $c_4=\frac{1}{AB}>0$.   
	\end{itemize}
	Therefore, if $c_0=c_1=0$, we have $\displaystyle \lim_{s\rightarrow 0^+} g(s)=A=0$.
	
	\underline{{(2) $\Leftrightarrow$ (4)}} We will see that  $\displaystyle \lim_{s\rightarrow 0^+} g(s)=0$ is equivalent to the induced metric being Bochner-K\"ahler. ordinary differential equation for a $U(n)$ invariant Bochner-K\"ahler metric is given by $sg'(s)=F(g(s))=c_4 g^3 +c_3 g^2 +g$ \cite{TL70}. Then we have $F(0)=0$, which implies $\lim\limits_{s\rightarrow 0 } g(s) = 0$. Converse is obvious.
	
\end{proof}

\begin{proof}[Proof of Corollary \ref{thm1}]
	The extremal equation is given by $sg'(s)=F(g(s))$ where $F(g)=c_4g^3+c_3 g^2 + g$.
	
	When $c_4=0$, the metric is cscK, and it follows from Theorem \ref{n smooth solution} that we cannot have $R<0$.
	
	Lemma \ref{lemma6.2} implies that, if we have $c_4\neq 0 $, then $\displaystyle\lim_{s\rightarrow +\infty} g(s) = B<\infty$ for degree reasons.
	
	The scalar curvature $R(s)$ is given by 
	$$R=-(n+2)(n+1)c_4 g(s) -(n+1)nc_3.$$
	The condition $R<0$ gives $-n c_3\leq (n+2)c_4 g(s)$. Let us check \eqref{metric beta}--\eqref{metric betagamma2} of Theorem \ref{Thm2.2} to see this is impossible.
	
	\begin{enumerate}[label=\underline{\smash{\textit{Case}\,(\arabic*)}}, align=left, leftmargin=0pt] \setcounter{enumi}{1}
		\item  We have $F(g) = c_4g(g-\beta)^2$ where $\beta=\displaystyle\lim_{s\rightarrow +\infty} g(s)$, $c_4 =\frac{1}{\beta^2}$, and $c_3=-\frac{2}{\beta}$. Then, $R<0$ implies $n\frac{2}{\beta}\leq (n+2)\frac{1}{\beta^2}g(s)$, which contradicts $\displaystyle\lim_{s\rightarrow 0^+} g(s) =0$.
		
		\item We have $F(g)=c_4 g(g-\beta)(g-\gamma)$ where $\gamma<0<\beta$, $\displaystyle\lim_{s\rightarrow +\infty}g(s) =\beta$, $c_4=\frac{1}{\beta\gamma}$, and $c_3= -\frac{\beta+\gamma}{\beta\gamma}$. Inequality $R<0$ implies $n\frac{\beta+\gamma}{\beta\gamma}\leq (n+2)\frac{1}{\beta\gamma}g(s)$. Since $\beta\gamma<0$, we have $g(s)\leq \frac{n+2}{n}g(s)\leq \beta+\gamma$. This contradicts $\displaystyle\lim_{s\rightarrow +\infty}g(s)=\beta$.
		
		\item We have $F(g) = c_4 g (g-\beta)(g-\gamma)$ where $0<\beta<\gamma$, $\displaystyle\lim_{s\rightarrow +\infty}g(s) =\beta$, $c_4 = \frac{1}{\beta\gamma}$, and $c_3=-\frac{\beta+\gamma}{\beta\gamma}$. Inequalities $R<0$ and $\beta\gamma>0$ imply $\frac{n}{n+2}(\beta+\gamma)\leq g(s)$. This contradicts $\displaystyle\lim_{s\rightarrow 0^+}g(s) =0$.
	\end{enumerate}
\end{proof}

\section{$U(2)$ Invariant Extremal K\"ahler Metrics on $\mathbb{C}^2 \backslash \{0\}$}\label{section 4}

\subsection{List of Solutions on $\mathbb{C}^2\backslash\{0\}$}

In \cite{HeLi}, He and Li give the list of scalar-flat and positive constant-scalar-curvature K\"ahler metrics defined on $\mathbb C^2 \backslash\{0\}$.  Here we extend their approach to \emph{extremal} K\"ahler metrics with non-constant scalar curvature on $\mathbb C^2 \backslash\{0\}$. Together with He and Li's Theorems \ref{2 zero singular solution}, \ref{2 negative singular solution}, our Theorem \ref{lemma2} classify solutions on $\mathbb{C}^2\backslash\{0\}$ to the $U(n)$ invariant extremal K\"ahler equation $sg'(s)=F(g(s))$, based on zeros of $F$.

We start this section with the theorem of He and Li, which gives the list of $U(2)$ invariant scalar-flat K\"ahler metrics on $\mathbb{C}^2\backslash\{0\}$:

\begin{theorem}[He-Li \cite{HeLi}, Theorem 1.2 ]\label{2 zero singular solution}
	Let $u:(0,+\infty) \rightarrow \mathbb{R}$ be a smooth function such that $u(s)$ is the potential of a K\"ahler metric with constant scalar curvature $R=0$ on $\mathbb{C}^2 \backslash \{0\}$. Then one of the following is true:
	
	\begin{enumerate}
		\item \label{Thm3.1.1} There exist constants $a,c$ with $a>0$ such that $g(s)=su'(s)$ is the smooth strictly increasing function $g:(0,\infty)\rightarrow (0,\infty)$ determined by the ordinary differential equation $sg'(s) = F(g(s))=g(s)$. Integrating, we get $g(s)=as$, hence $u(s)=as+c$.
		
		\item \label{Thm3.1.2} There exist constants $\alpha,a,c$ with $\alpha>0$, $a>0$, such that $g(s)=su'(s)$ is the smooth strictly increasing function $g:(0,\infty)\rightarrow(\alpha,\infty)$ determined by the ordinary differential equation $sg'(s) = F(g(s))=g(s)-\alpha$. Integrating, we get $g(s)=\alpha + as$, hence $u(s) = \alpha\log s + as + c$.
		
		\item \label{Thm3.1.3} There exist constants $\alpha,\beta,c$ with $\alpha \neq 0, \beta>0, \alpha < \beta$
		such that $g(s)=su'(s)$ is the smooth strictly increasing function $g:(0,\infty)\rightarrow(\beta,\infty)$ determined by the ordinary differential equation $sg'(s) = F(g(s)) = \frac{(g(s)-\alpha)(g(s)-\beta)}{g(s)}$. Integrating, we get
		\[
		\frac{\beta}{\beta-\alpha} \log(g(s)-\beta) -\frac{\alpha}{\beta-\alpha} \log(g(s)-\alpha)=\log s +c.
		\]
		
		\item \label{Thm3.1.4} There exist constants $\alpha,c$ with $\alpha >0$ such that $g(s)=su'(s)$ is the smooth strictly increasing function $g: (0, \infty)\rightarrow (\alpha,\infty)$ determined by the ordinary differential equation $sg'(s) = F(g(s)) = \frac{(g(s)-\alpha)^2}{g(s)}$. Integrating, we get
		\[
		\log(g(s)-\alpha)-\frac{\alpha}{g(s)-\alpha}=\log s +c.
		\]
	\end{enumerate}
\end{theorem}

The next theorem\footnote{In \cite{HeLi}, the authors state this theorem for a particular value of the scalar curvature, namely $R=6$. Here, we present the statement for a general constant $R>0$.} gives the list of  $U(2)$ invariant  positive cscK metrics on $\mathbb{C}^2 \backslash \{0\}$. 

\begin{theorem}[He-Li \cite{HeLi}, Theorem 1.3]\label{2 negative singular solution}
	Let $u:(0,+\infty) \rightarrow \mathbb{R}$ be a smooth function such that $u(s)$ is the potential of a K\"ahler metric with positive constant scalar curvature $R$ on $\mathbb{C}^2\backslash\{0\}$ and $g(s)=su'(s)$. Then one of the following is true:
	
	\begin{enumerate}
		\item \label{Thm3.2.1} There exist constants $a,c$ with $a>0$ such that
		\[
		u(s)=\log(s+a)+c
		\]

		\item \label{Thm3.2.2}There exist constants $\alpha,\beta, c$ with $0<\alpha<\beta<\infty$ such that $g(s)$ is the smooth strictly increasing function $g:(0,\infty)\rightarrow (\alpha,\beta)$ determined by the ordinary differential equation $sg'(s)=F(g(s))=-\frac{1}{\alpha+\beta} (g(s)-\alpha)(g(s)-\beta)$. Integrating, we get
		\[
		\log(g(s)-\alpha)-\log(\beta-g(s)) = \frac{\beta-\alpha}{\alpha+\beta} \log s + c.
		\]
		
		\item \label{Thm3.2.3} There exist constants $\alpha,\beta,\gamma,c$ with 
		$\alpha \neq 0, \beta>0, \alpha<\beta<\gamma, \alpha+\beta+\gamma>0$ 
		such that $g(s)$ is the smooth strictly increasing function $g:(0,\infty)\rightarrow(\beta,\gamma)$ determined by the ordinary differential equation $sg'(s) = F(g(s)) = -\frac{1}{\alpha+\beta+\gamma} \frac{(g(s)-\alpha)(g(s)-\beta)(g(s)-\gamma)}{g(s)}$. Integrating, we get
		
		\begin{equation*}
		\begin{split}
		&-\alpha(\gamma-\beta)\log(g(s)-\alpha)+\beta(\gamma-\alpha)\log(g(s)-\beta) \\
		&-\gamma(\beta-\alpha)\log(\gamma-g(s)) 
		=\frac{(\beta-\alpha)(\gamma-\beta)(\gamma-\alpha)}{\alpha+\beta+\gamma}\log s+c.
		\end{split}	
		\end{equation*}	
		In this case, we have $R=\frac{6}{\alpha+\beta+\gamma}$.
		
		\item \label{Thm3.2.4} There exist  constants $\alpha,\beta,c$ with $0<\alpha<\beta$ 
		such that $g(s)$ is the smooth strictly increasing function $g:(0,\infty)\rightarrow(\alpha,\beta)$ determined by the ordinary differential equation $sg'(s) = F(g(s)) = -\frac{1}{2\alpha+\beta} \frac{(g(s)-\alpha)^2(g(s)-\beta)}{g(s)}$. Integrating, we get
		\begin{equation*}
		-\frac{\alpha(2\alpha+\beta)}{\beta - \alpha} \frac{1}{g-\alpha} + \frac{\beta (2\alpha + \beta)}{(\beta-\alpha)^2} \log(g-\alpha) - \frac{\beta(2\alpha+\beta)}{(\beta-\alpha)^2} \log(\beta - g) = \log s + c.
		\end{equation*}
		In this case, $R=\frac{6}{2\alpha+\beta}$.
	\end{enumerate}
\end{theorem}

In the next theorem, we expand the above list of metrics to contain  all $U(2)$ invariant extremal K\"ahler metrics on $\mathbb{C}^2\backslash\{0\}$.

\begin{theorem}\label{lemma2}
	Let $u:(0,\infty)\rightarrow \mathbb{R}$ be a smooth function such that $u(s)$ is the potential of a K\"ahler metric satisfying the extremal condition on $\mathbb{C}^2\backslash \{0\}$ and $g(s)=su'(s)$. Then one of the following is true.
	\begin{enumerate}
		\item Metric can be extended smoothly to $\mathbb{C}^2$.
		\item Metric is cscK with a singularity at the origin.
		\item There exist constants $\alpha,\beta,c$ with $0<\alpha<\beta$ such that $g(s)$ is the smooth strictly increasing function $g:(0,\infty)\rightarrow(\alpha,\beta)$ determined by $sg'=F(g)=\frac{1}{\beta(\beta+2\alpha)} (g-\alpha)(g-\beta)^2$. Integrating, we get
		\begin{align}\label{*11}
		\frac{\beta(\beta+2\alpha)}{(\alpha-\beta)^2} \left\{\log(g(s)-\alpha) - \log(\beta-g(s)) - \frac{\beta-\alpha}{g(s)-\beta}\right\} = \log s + c.
		\end{align}
		
		\item \label{case4} There exist constants $\alpha,\beta,\gamma, c$ with $\gamma <-\frac{\alpha\beta}{\alpha+\beta} <0<\alpha<\beta$ such that $g(s)$ is the smooth strictly increasing function $g:(0,\infty)\rightarrow(\alpha,\beta)$ determined by $sg'=F(g)=\frac{1}{\alpha\beta +\alpha\gamma +\beta\gamma} (g-\alpha)(g-\beta)(g-\gamma)$. Integrating, we get
		\begin{align} \label{*10}
		(\alpha\beta +\alpha\gamma +\beta\gamma) \left\{\frac{\log(g(s)-\alpha)}{(\alpha-\beta)(\alpha -\gamma)} + \frac{\log(\beta -g(s))}{(\beta-\alpha)(\beta -\gamma)} \right.\\ \left.\nonumber + \frac{\log|g(s)-\gamma|}{(\gamma-\alpha)(\gamma-\beta)} \right\}=\log s + c.
		\end{align}
		
		\item There exist constants $\alpha,\beta,\gamma, c$ with $0<\alpha<\beta<\gamma$ such that  $g(s)$ is the smooth strictly increasing function $g:(0,\infty)\rightarrow(\alpha,\beta)$ determined by equation $sg'=F(g)$, where $F$ is given as in Case \eqref{case4}. Integrating, we get \eqref{*10}.
		
		\item \label{case6}There exist constants $\alpha,\beta, c$ with $0<\alpha<\beta$ such that $g(s)$ is the smooth strictly increasing function $g:(0,\infty)\rightarrow(\alpha,\beta)$  determined by  $sg'=F(g)=\frac{1}{(\alpha+\beta)^2 + 2\alpha\beta} \frac{ (g-\alpha)^2(g-\beta)^2}{g}$. Integrating, we get
		\begin{align}\label{*6}
		((\alpha+\beta)^2 + 2\alpha\beta) \left\{ \frac{\alpha+\beta}{(\beta-\alpha)^3} \log(g(s)-\alpha) - \frac{\alpha}{(\alpha-\beta)^2} \frac{1}{g(s)-\alpha} - \right.\\ \nonumber \left. \frac{\alpha+\beta}{(\beta-\alpha)^3} \log(\beta - g(s)) - \frac{\beta}{(\alpha-\beta)^2} \frac{1}{g(s)-\beta} \right\} = \log s + c.
		\end{align}
		
		\item \label{case7} There exist constants $\alpha,\beta,\gamma, c$ with $0<\alpha<\beta<\gamma$ such that  $g(s)$ is the smooth strictly increasing function $g:(0,\infty)\rightarrow(\alpha,\beta)$  determined by  $sg'=F(g)=\frac{1}{\alpha^2 +2\alpha\beta+2\alpha\gamma +\beta\gamma} \frac{ (g-\alpha)^2(g-\beta)(g-\gamma)}{g}$. Integrating, we get
		\begin{align}\label{*3}
		(\alpha^2 +2\alpha\beta+2\alpha\gamma +\beta\gamma) \left\{ \frac{\gamma}{(\gamma-\alpha)(\gamma-\beta)} \log(\gamma-g(s)) - \right. \\\nonumber \left. \frac{\alpha}{(\alpha-\gamma)(\alpha-\beta)} \frac{1}{g(s)-\alpha} + \frac{-\alpha^2+\beta\gamma}{(\alpha-\gamma)(\alpha-\beta)} \log(g(s)-\alpha) +\right. \\\nonumber \left. \frac{\beta}{(\beta-\gamma)(\beta-\alpha)^2} \log|\beta-g(s)| 
		\right\} = \log s + c.
		\end{align}
		
		\item There exist constants $\alpha,\beta,\gamma, c$ with $\alpha < - \frac{\alpha^2 +\beta\gamma}{2(\beta+\gamma)} <0<\beta<\gamma$ such that  $g(s)$ is the smooth strictly increasing function $g:(0,\infty)\rightarrow(\beta,\gamma)$ determined by equation $sg'=F(g)$, where $F(g)$ is given as in Case \eqref{case7}. Integrating, we get \eqref{*3}.
		
		\item There exist constants $\alpha,\beta,\gamma, c$ with $\alpha < - \frac{\beta^2 +2\beta\gamma}{2\beta+\gamma} <0<\beta<\gamma$ such that  $g(s)$ is the smooth strictly increasing function $g:(0,\infty)\rightarrow(\beta,\gamma)$ determined by $sg'=F(g)=\frac{1}{\beta^2 +2\beta\alpha+2\beta\gamma +\alpha\gamma} \frac{ (g-\alpha)(g-\beta)^2(g-\gamma)}{g}$. Integrating, we get
		\begin{align}\label{*4}
		(\beta^2 +2\beta\alpha+2\beta\gamma +\alpha\gamma) \left\{ \frac{\alpha}{(\alpha-\beta)(\alpha -\gamma)} \log(g(s)-\alpha) \right.\\ \nonumber \left. - \frac{\beta}{(\beta-\alpha)(\beta-\gamma)} \frac{1}{g(s)-\beta} + \frac{-\beta^2+\alpha\gamma}{(\beta-\alpha)(\beta-\gamma)} \log(g(s)-\beta) \right. \\\nonumber \left. + \frac{\gamma}{(\gamma-\alpha)(\gamma-\beta)^2} \log(\gamma-g(s)) 
		\right\} = \log s + c.
		\end{align}
		
		\item \label{case10} There exist constants $\alpha,\beta,\gamma, c$ with $-\frac{\gamma(\gamma+2\beta)}{2\gamma+\beta}<\alpha < \beta<\gamma$ and $\alpha\beta\gamma\neq 0$ such that  $g(s)$ is the smooth strictly increasing function $g:(0,\infty)\rightarrow(\beta,\gamma)$ determined by  $sg'=F(g)=\frac{1}{\gamma^2 +2\beta\gamma+2\alpha\gamma +\alpha\beta} \frac{ (g-\alpha)(g-\beta)(g-\gamma)^2}{g}$. Integrating, we get
		\begin{align}\label{*5}
		(\gamma^2 +2\beta\gamma+2\alpha\gamma +\alpha\beta) \left\{ \frac{\alpha}{(\alpha-\beta)(\alpha -\gamma)} \log(g(s)-\alpha) \right.\\ \nonumber \left. - \frac{\gamma}{(\gamma-\alpha)(\gamma-\beta)} \frac{1}{g(s)-\gamma} + \frac{-\gamma^2+\alpha\beta}{(\gamma-\alpha)(\gamma-\beta)} \log(\gamma-g(s)) \right. \\\nonumber \left. + \frac{\beta}{(\beta-\alpha)(\beta-\gamma)^2} \log(g(s)-\beta) 
		\right\} = \log s + c.
		\end{align}
		
		\item \label{case11}There exist constants $\alpha,\beta,\gamma, \tau, c$ with $-\frac{\beta\gamma+\beta\tau+\gamma\tau}{\beta+\gamma+\tau}<\alpha < \beta<\gamma<\tau$ and $\alpha\beta\gamma\tau\neq 0$ such that  $g(s)$ is the smooth strictly increasing function $g:(0,\infty)\rightarrow(\beta,\gamma)$ determined by $$sg'=F(g)=\frac{1}{\alpha\beta+\alpha\gamma+\alpha\tau+\beta\gamma+\beta\tau+\gamma\tau} \frac{ (g-\alpha)(g-\beta)(g-\gamma)(g-\tau)}{g}.$$ Integrating, we get
		\begin{align}\label{**6}
		(\alpha\beta+\alpha\gamma+\alpha\tau+\beta\gamma+\beta\tau+\gamma\tau) \left\{ \frac{\alpha}{(\alpha-\beta)(\alpha-\gamma)(\alpha-\tau)}\log(g(s)-\alpha)  \right. \\\nonumber \left. +
		\frac{\beta}{(\beta-\alpha)(\beta-\gamma)(\beta-\tau)}\log|g(s)-\beta|+
		\frac{\gamma}{(\gamma-\alpha)(\gamma-\beta)(\gamma-\tau)}\log|g(s)-\gamma| \right. \\\nonumber \left. +
		\frac{\tau}{(\tau-\alpha)(\tau-\beta)(\tau-\gamma)}\log|g(s)-\tau| 
		\right\}
		=\log s + c .
		\end{align}
		
		\item There exist constants $\alpha,\beta,\gamma, \tau, c$ with $\alpha\beta+\alpha\gamma+\alpha\tau+\beta\gamma+\beta\tau+\gamma\tau<0$, $\alpha\beta\gamma\tau\neq 0$ such that  $g(s)$ is the smooth strictly increasing function $g:(0,\infty)\rightarrow(\gamma,\tau)$ determined by $sg'=F(g)$ where $F$ is in Case \eqref{case11}. Integrating, we get \eqref{**6}.
		
		\item There exist constants $\alpha,\beta,a,b$ with $0<\alpha<\beta$ and $a^2+2a(\alpha+\beta)+b^2+\alpha\beta<0$ such that  $g(s)$ is the smooth strictly increasing function $g:(0,\infty)\rightarrow(\alpha,\beta)$ determined by $$sg'=F(g)=\frac{1}{a^2 +2a(\alpha+\beta)+b^2+\alpha\beta} \frac{ (g-\alpha)(g-\beta)(g-2ag+(a^2+b^2))}{g}.$$ Integrating, we get
		\begin{align}\label{**7}
		(a^2 +2a(\alpha+\beta)+b^2+\alpha\beta) \Bigl\{c_1 \log(g(s)-\alpha) +c_2 \log(\beta -g(s)) +  \Bigr. \\\nonumber \Bigl. \int_{\alpha}^{\beta}\frac{-(c_1 + c_2) g(s) + 2(c_1+c_2) a - c_1 \alpha -c_2 \beta}{g^2(s) -2a g(s)+a^2+b^2} g'(s) ds\Bigr\} = \log s,
		\end{align} 
		where $c_1 = \frac{\alpha}{(\alpha-\beta)(\alpha^2-2a\alpha +a^2+b^2)}$ and $c_2 = \frac{\beta}{(\beta -\alpha)(\beta^2-2a\beta +a^2+b^2)}$. 
	\end{enumerate}
	
\end{theorem}

\begin{proof}
	See Section \ref{section proof 3}.
\end{proof}

\subsection{Examples of Extremal K\"ahler Metrics on Line Bundles over $\mathbb{CP}^1$}

The list of $U(2)$ invariant extremal K\"ahler metrics on $\mathbb{C}^2\backslash\{0\}$ given by Theorems \ref{2 zero singular solution}, \ref{2 negative singular solution}, and \ref{lemma2} can be used to write down  extremal K\"ahler metrics on $\mathbb{CP}^2$, Hirzebruch surfaces $\mathcal{F}_k^2$, $k\geq 1$, and line bundles over $\mathbb{CP}^1$.

On $\mathcal{O}(-k)$ bundles, examples of $U(n)$ invariant scalar-flat K\"ahler metrics were given by Burns, Eguchi-Hanson, LeBrun, Pedersen-Poon, and Simanca. (See Remark \ref{remark 3.4} and Introduction for details). Apostolov-Rollin generalized these metrics to scalar-flat metrics on $\mathbb{CP}^n _{[-a_0,a_1,\dots,a_n]}$. In this section, we retrieve scalar-flat metrics on $\mathcal{O}_{\mathbb{CP}^1}(-k)$, $k\geq 1$, using solutions listed in Theorem \ref{2 zero singular solution}.

$U(2)$ invariant extremal K\"ahler metrics with positive constant scalar curvature on $\mathbb{C}^2\backslash\{0\}$ are listed in Theorem \ref{2 negative singular solution} \cite{HeLi}. For a list of smooth and singular metrics induced by Theorem \ref{2 negative singular solution} on $\mathbb{CP}^2$, Hirzebruch surfaces, and $\mathcal{O}_{\mathbb{CP}^1}(k)$, $k\geq 1$,  see Remark \ref{rmk-pos-csck}. Metrics with positive constant-scalar-curvature on $\mathcal{O}_{\mathbb{CP}^1}(k)$ seem to be new.

Some of the metrics induced by Theorem \ref{lemma2} include Calabi's extremal metrics with non-constant scalar curvature on Hirzebruch surfaces, complete strictly extremal metrics on $\mathcal{O}_{\mathbb{CP}^1}(-k)$, $k\geq 1$, and $\mathbb{C}^2\backslash\{0\}$. To the best of our knowledge the last two metrics are new.

\begin{remark}\label{remark 3.4}
	In \cite{LeBrun88}, LeBrun solves the equation for $U(2)$ invariant scalar-flat K\"ahler metrics on $\mathbb{C}^2\backslash\{0\}$. We can retrieve his results from Theorem \ref{2 zero singular solution} as follows. 
	\begin{enumerate}[label=(\textit{\roman*})]
		\item  Case \eqref{Thm3.1.1} induces a multiple of the Euclidean metric on $\mathbb{C}^2$.
		
		\item Case \eqref{Thm3.1.2} induces Burns metric on  $\mathcal{O}_{\mathbb{CP}^1}  (-1)$.
		
		\item Case \eqref{Thm3.1.3} induces LeBrun metrics on $\mathcal{O}_{\mathbb{CP}^1}(-k)$, ($k\geq 2$). When $k=2$, this gives the Ricci-flat Eguchi-Hanson metric on $\mathcal{O}_{\mathbb{CP}^1}(-2)$. 
		
		In this case, the solutions are determined by $F(g)=\frac{(g-\alpha) (g-\beta)}{g}$, where $\alpha\neq 0$, $0<\beta$, $\alpha<\beta$, $\beta<g(s)<\infty$. When $\alpha,\beta$ are such that $F'(\beta) = \frac{\beta-\alpha}{\beta} = k$ ($k\geq 2$) holds, solutions induce scalar-flat K\"ahler ALE metrics on $\mathcal{O}_{\mathbb{CP}^1}(-k)$.  When $k=2$,  we have $\beta=-\alpha$, which induces the Eguchi Hanson metric.
		
		\item Case \eqref{Thm3.1.4} induces  complete scaler-flat K\"ahler metrics on $\mathbb{C}^2\backslash\{0\}$. Completeness as $|z|\rightarrow 0$ and $|z|\rightarrow \infty$ can be checked as in Example \ref{ex4.6}.
	\end{enumerate}
\end{remark}

\begin{remark}\label{rmk-pos-csck}
	The solutions given by Theorem \ref{2 negative singular solution} induce positive cscK metrics on various complex surfaces. It is well known that there are no smooth cscK metrics on Hirzebruch surfaces $\mathcal{F}^2_k$, $k\geq 1$.  
	
	\begin{enumerate}[\itshape(i)]
	\item Case \eqref{Thm3.2.1} gives the Fubini-Study metric on $\mathbb{CP}^2$. 
	
	\item Case \eqref{Thm3.2.2} gives the cscK metrics on Hirzebruch surfaces $\mathcal{F}^2_k$ with cone angles $2\pi\theta/k$, $0<\theta<1$, along $0$- and $\infty$-sections. In this case we have $F(g) = -\frac{1}{\alpha+\beta} (g-\alpha)(g-\beta)$ such that $0<\alpha<g(s)<\beta$. Here we have $\theta=F'(\alpha)= -F'(\beta) = \frac{\beta-\alpha}{\alpha+\beta}$.
	
	\item Case \eqref{Thm3.2.4} gives complete positive cscK metrics on $\mathcal{O}_{\mathbb{CP}^1}(k)$, $k\geq 1$, which seems to be new, and we show it in the next example. 
	
	\end{enumerate}
\end{remark}

The next example gives positive cscK  metrics on  $\mathcal{O}_{\mathbb{CP}^1}(k)$, $k\geq 1$.

\begin{example}[Positive cscK metrics on $\mathcal O_{\mathbb{CP}^1}(k)$, $k\geq 1$] \label{ex4.6}	
		Let us consider the positive cscK metric on $\mathbb C^2\backslash\{0\}$ given by the ordinary differential equation
		\begin{equation}\label{dagger}
		sg'(s) = F(g(s)) = -\frac{1}{2\alpha+\beta}\frac{(g(s)-\alpha)^2(g(s)-\beta)}{g(s)},
		\end{equation}
		where we have $0<\alpha<\beta$. We have $\displaystyle\lim_{s\rightarrow 0^+} g(s)=\alpha$ and $\displaystyle\lim_{s\rightarrow +\infty} g(s) = \beta$. 
		Integrating, there exists a constant $c$ such that
		\begin{equation}\label{newstar}
		-\frac{\alpha(2\alpha+\beta)}{\beta - \alpha} \frac{1}{g-\alpha} + \frac{\beta (2\alpha + \beta)}{(\beta-\alpha)^2} \log(g-\alpha) - \frac{\beta(2\alpha+\beta)}{(\beta-\alpha)^2} \log(\beta - g) = \log s + c.
		\end{equation}
		
		 We can obtain the positive line bundle $\mathcal{O}_{\mathbb{CP}^1}(k)$, $k>0$, by gluing a $\mathbb{CP}^1$ to $(\mathbb{C}^2\backslash\{0\})\slash \mathbb{Z}_k$ at $|z|=\infty$. $U(2)$ invariant  K\"ahler metrics  on $\mathbb{C}^2\backslash\{0\}$ induce metrics on $\mathcal{O}_{\mathbb{CP}^1}(k)\backslash S_{\infty}$. Here $S_\infty$ stands for the zero section of $\mathcal{O}_{\mathbb{CP}^1}(k)$. The induced metric can be extended by continuity to a smooth metric on $\mathcal{O}_{\mathbb{CP}^1}(k)$ if and only if $F(\beta)=0$ and $F'(\beta)= -k$.
		
		The condition $F(\beta)=0$ is clearly satisfied. We compute $F'(\beta)=-\frac{(\beta-\alpha)^2}{\beta^2 (\beta+2\alpha)}$. We note that for every positive integer $k$ there exist constants $\alpha,\beta$, $0<\alpha<\beta$ which satisfy 
		$$\frac{(\beta-\alpha)^2}{\beta^2 (\beta+2\alpha)}=k.$$
		In order to see that such $\alpha,\beta$ exist for each $k\geq 1$, let us introduce new variables $x=\alpha>0$ and $y=\beta-\alpha>0$. Then the above equation becomes 
		$$y^2-k(x+y)^2 (y+3x) = 0.$$
		For each positive integer $k$, this equation has solutions $(x,y)$ with $x>0$, $y>0$. 
		
		The function $F(g)$ is strictly positive on $(\alpha,\beta)$. It follows from Calabi \cite{Calabi82} that the K\"ahler metric extends smoothly to $\mathcal{O}_{\mathbb{CP}^1}(k)$.
		
		We need to show that the induced metric is complete on the total space of $\mathcal{O}_{\mathbb{CP}^1}(k)$.
		
		The metric is complete if the improper integral that gives the geodesic distance to $z=0$ 	
		$$\int_{0}^{s_0} \sqrt{\frac{g'(s)}{s}} ds = \int_0^{s_0} \frac{\sqrt{F(g(s))}}{s} ds$$
		is infinite.
		
		Since $\displaystyle\lim_{s\rightarrow 0^+}g(s)=\alpha>0$ and $\displaystyle\lim_{s\rightarrow \infty}g(s)=\beta$, $\frac{\beta-g}{(2\alpha+\beta)g}$ is bounded on $(0,s_0)$. There exists $d_1>0$ such that 
		$$\sqrt{F(g)} = \sqrt{\left(-\frac{1}{2\alpha+\beta}\right) (g-\alpha)^2 (g-\beta) \frac{1}{g}} >d_1 (g-\alpha).$$
		Then, we have
		$$\int_0^{s_0} \sqrt{\frac{g'(s)}{s} } ds>d_1 \int_{0}^{s_0} \frac{g-\alpha}{s} ds.$$
		If we choose $s_0$ small enough, we have $\log(g-\alpha)<0$ on $(0,s_0)$, and $\log(\beta-g)$ is bounded. Therefore, Equation \eqref{newstar} implies 
		\begin{align*}
		-\frac{\alpha(2\alpha+\beta)}{\beta-\alpha}\frac{1}{g(s)-\alpha} > \log s + c\\
		\frac{g(s)-\alpha}{s}> -\frac{\beta-\alpha}{\alpha(2\alpha+\beta)} \frac{1}{s (\log s+ c)}.
		\end{align*}
		Integrating both sides of this inequality on $(0,s_0)$ we see that the integral
		$$ \int_{0}^{s_0} \frac{g(s)-\alpha}{s} ds$$
		is infinite.
\end{example}

The next example shows that Case \ref{case10} of Theorem \ref{lemma2} induces an extremal  K\"ahler metric on  $\mathcal O_{\mathbb{CP}^1} (-k)$, $k\geq 1$ with non-constant scalar curvature.

\begin{example}[Strictly extremal metrics on $\mathcal O_{\mathbb{CP}^1} (-k)$, $k\geq 1$]\label{example 3.7}
		Let us consider Case \ref{case10} of Theorem \ref{lemma2}. Since $c_4\neq 0$, it follows from Equation \eqref{R} that this is a strictly extremal metric on $\mathbb{C}^2\backslash \{0\}$. Extremal equation is given by $sg'(s)=F(g(s))$ where
		$$F(g)=\frac{c_4(g-\alpha)(g-\beta)(g-\gamma)^2}{g}.$$
		Here, we have $\alpha<\beta<\gamma$, $\alpha\beta\gamma\neq 0$, $\displaystyle
		\lim_{s\rightarrow 0^+}g(s)=\beta$, $\displaystyle
		\lim_{s\rightarrow \infty}g(s)=\gamma$ and $c_4 = \displaystyle\frac{1}{\alpha\beta +2\alpha\gamma+2\beta\gamma+\gamma^2}>0$. We note that $c_4>0$ implies $-\frac{\gamma(\gamma+2\beta)}{2\gamma + \beta} <\alpha$.
		
		As in the proof of Theorem \ref{lemma2}, we will rewrite the ordinary differential equation as
		\begin{align*}
		g'(\gamma^2 +2\beta\gamma+2\alpha\gamma +\alpha\beta) \left\{ \frac{\alpha}{(\alpha-\beta)(\alpha -\gamma)} \frac{1}{g(s)-\alpha} \right.\\ \nonumber \left. + \frac{\gamma}{(\gamma-\alpha)(\gamma-\beta)} \frac{1}{(g(s)-\gamma)^2} + \frac{-\gamma^2+\alpha\beta}{(\gamma-\alpha)(\gamma-\beta)} \frac{1}{g(s)-\gamma} \right. \\\nonumber \left. + \frac{\beta}{(\beta-\alpha)(\beta-\gamma)^2} \frac{1}{g(s)-\beta}
		\right\} = \frac{1}{s}.
		\end{align*}
		Now recall that we can obtain the line bundle $\mathcal{O}_{\mathbb{CP}^1}(-k)$, $k=1,2,\dots$, by gluing a $\mathbb{CP}^1$ to $(\mathbb{C}^2\backslash \{0\})\slash \mathbb{Z}_k$ at $z=0$. The $U(2)$ invariant K\"ahler metric on $\mathbb{C}^2\backslash \{0\}$ determined by $sg'(s)=F(g(s))$ induces a metric on $\mathcal{O}_{\mathbb{CP}^1}(-k)\backslash S_0$. The induced metric can be extended by continuity to a smooth metric on $\mathcal{O}_{\mathbb{CP}^1}(-k)$ if and only if $F(\beta)=0$ and $F'(\beta)=k$. The condition $F(\beta)=0$ is clearly satisfied. We compute
		$$F'(\beta)= \frac{(\beta-\alpha)(\gamma-\beta)^2}{\beta(\alpha\beta +2\alpha\gamma+2\beta\gamma+\gamma^2)}.$$
		We need to show that for every positive integer $k$, there exist constants $\alpha,\beta,\gamma$ with 
		\begin{equation}\label{*}
		-\frac{\gamma(\gamma+2\beta)}{2\gamma + \beta} <\alpha<\beta<\gamma, \quad \alpha\beta\gamma\neq 0, 
		\end{equation}
		that satisfy
		\begin{equation}\label{**}
		F'(\beta) = \frac{(\beta-\alpha)(\gamma-\beta)^2}{\beta(\alpha\beta +2\alpha\gamma+2\beta\gamma+\gamma^2)} =k.
		\end{equation}
		For simplicity, let $\gamma=2\beta$. Then, Equations \eqref{*} and \eqref{**} give
		\begin{equation}\label{Delta}
		-\frac{8\beta}{5}<\alpha<\beta, \quad \alpha\beta\neq 0, \quad \frac{\beta-\alpha}{5\alpha+8\beta} =k
		\end{equation}
		For each positive integer $k$, the pair $(\alpha,\beta)=\left(\frac{1-8k}{1+5k}\beta,\beta\right)$ satisfies \eqref{Delta}.
		
		We need to show that the induced metric on $\mathcal{O}_{\mathbb{CP}^1}(-k)$ is complete as $|z|\rightarrow \infty$, i.e. as $g(s)\rightarrow\gamma$. The metric is complete if the improper integral
		$$\int_{s_0}^{\infty} \sqrt{\frac{g'(s)}{s}}\ ds = \int_{s_0}^{\infty} \frac{\sqrt{F(g(s))}}{s}\ ds $$
		is infinite. Since $\displaystyle\lim_{s\rightarrow 0^+} g(s)=\beta$ and $\displaystyle\lim_{s\rightarrow +\infty} g(s)=\gamma>0$, $\frac{c_4(g-\alpha)(g-\beta)}{g}$ is bounded on $(s_0,\infty)$. There exists $d_1>0$ such that $\sqrt{F(g)}>d_1 (\gamma-g)$. Then we have
		$$\int_{s_0}^{\infty} \sqrt{\frac{g'(s)}{s}}\ ds > d_1 \int_{s_0}^{\infty} \frac{\gamma-g(s)}{s} \ ds. $$
		If we choose $s_0$ large enough, we have $\log(\gamma-g)<0$ on $(s_0,\infty)$, and $\log(g-\alpha)$, $\log(g-\beta)$ are bounded. Noting that $-\gamma^2+\alpha\beta<0$, Equation \eqref{*5} implies 
		\begin{align*}
		\frac{\gamma}{(\gamma-\alpha)(\gamma-\beta) }\frac{1}{\gamma-g(s)} &< \frac{1}{\gamma^2+2\beta\gamma + 2\alpha\gamma + \alpha\beta} \log s + c_1 \\
		\frac{\gamma-g(s)}{s} &> \frac{c_2\log s + c_3}{s},\quad c_2>0.
		\end{align*}
		Integrating both sides of this inequality on $(s_0,\infty)$, we see that the integral
		$$\int_{s_0}^{\infty}\frac{\gamma-g(s)}{s}\ ds$$
		is infinite.
\end{example}

As our last example, we give a complete extremal K\"ahler metric on $\mathbb{C}^2\backslash \{0\}$.

\begin{example}[Strictly extremal metrics on $\mathbb{C}^2\backslash \{0\}$] \label{ex4.8}
	Let us consider Case \eqref{case6} of Theorem \ref{lemma2}. Since $c_4\neq 0$, it follows from Equation \eqref{R} that this is a strictly extremal metric on $\mathbb{C}^2\backslash \{0\}$. Extremal equation is given by $sg'(s)=F(g(s))$ where
	$$F(g)=\frac{1}{(\alpha+\beta)^2+2\alpha\beta} \frac{(g-\alpha)^2(g-\beta)^2}{g}.$$
	Here, we have $0<\alpha<\beta$,  $\displaystyle
	\lim_{s\rightarrow 0^+}g(s)=\alpha$, $\displaystyle
	\lim_{s\rightarrow \infty}g(s)=\beta$. It is straightforward computation as in  previous examples to check that the metric is complete as  $|z|\rightarrow 0 $ and $|z|\rightarrow \infty$. This gives us a complete extremal K\"ahler metric with non-constant scalar curvature on $\mathbb{C}^2 \backslash \{0\}$.
\end{example}

\subsection{Proofs}\label{section proof 3}

\begin{proof}[Proof of Theorem \ref{lemma2}]
	It follows from equation \eqref{ordinary differential equation} that there exists constants $c_0,c_1,c_3,c_4$ such that 
	\begin{align}\label{*8}
	\frac{gg'}{c_4 g^4 + c_3 g^3 + g^2 +c_1 g +c_0}= \frac{1}{s}.
	\end{align}
	
	\begin{enumerate}[label=\underline{\smash{\textit{Case}\,(\arabic*)}}, align=left, leftmargin=0pt]
		
		\item $c_4=0$.\\ 
		We see from equation \eqref{R} that $\omega$ is a cscK metric. Classification of cscK metrics on $\mathbb{C}^2\backslash 0$ is given by Theorem 1.2 in \cite{HeLi}.\\
		
		\item $c_4\neq 0$. $c_0 =c_1=0$ and $\lim_{s\rightarrow 0^+} g(s) = 0$.\\
		It follows from Remark that, in this case, the metric can be smoothly extended to the origin, hence Theorem \ref{Thm2.2} applies.\\
		
		\item $c_4\neq 0$, $c_0=c_1=0$ and $\lim_{s\rightarrow 0^+} g(s)=A>0$.\\
		We have $H(x) = c_4 x^3 +c_3 x^2 +x$ and it follows from Lemma \ref{lemma6.2} that $H(A)=0$, $B<\infty$, and $H(x)>0$ in $(A,B)$. In this case, the roots are given by $\gamma=0<\alpha =A <\beta=B$. But $H(x)>0$ on $(\alpha,\beta)$, and this implies that $c_4 = \frac{1}{\alpha\beta}<0$, which is a contradiction.
		
		\item $c_4\neq 0$, $c_0=0$, $c_1\neq 0$, and the polynomial $c_4x^3 +c_3 x^2 +x+c_1$ has roots $\alpha,\alpha,\alpha$.\\
		It follows from Lemma \ref{lemma6.2} that $B<\infty$ for degree reasons, and this case is impossible.\\
		
		\item $c_4\neq 0$, $c_0 =0$, $c_1\neq 0$, and the polynomial $c_4x^3 +c_3 x^2 +x+c_1$ has roots $\alpha,\alpha,\beta$ with $\alpha<\beta$.\\
		It follows from Lemma \ref{lemma6.2} that $B<\infty$, $\alpha = A>0$, $\beta =B$, and $H(x)>0$ on $(\alpha,\beta)$, which implies that $c_4<0$.
		
		Since $H(x)= c_4 (x-\alpha)^2 (x-\beta)$, we have $1=c_4 (\alpha^2 +2\alpha\beta)$, which contradicts to $0<\alpha<\beta$.\\
		
		\item  $c_4\neq 0$, $c_0 =0$, $c_1\neq 0$, and the polynomial $c_4x^3 +c_3 x^2 +x+c_1$ has roots $\alpha,\beta,\beta$ with $\alpha<\beta$.\\
		We have $\alpha\neq 0$, $\beta\neq 0$. The equation \eqref{*8} can be written as
		\begin{align*}
		\frac{g'(s)}{c_4 (\alpha-\beta)^2} \left\{\frac{1}{g(s)-\alpha} -\frac{1}{g(s)-\beta}  + \frac{\beta-\alpha}{(g(s)-\beta)^2}\right\} = \frac{1}{s}.
		\end{align*}
		It follows from Lemma \ref{lemma6.2} that $\alpha =A$ $\beta=B$, and $H(x)>0$ on $(\alpha,\beta)$. Hence $c_4=\frac{1}{\beta(\beta+2\alpha)}>0$, and there exists a constant $c$ such that 
		\begin{align}\tag{\ref{*11}}
		\frac{\beta(\beta+2\alpha)}{(\alpha-\beta)^2} \left\{\log(g(s)-\alpha) - \log(\beta -g(s)) -\frac{\beta-\alpha}{g(s)-\beta} \right\}= \log s + c.
		\end{align}
		On the other hand, Lemma \ref{lemma6.1} implies that there exists a unique smooth, strictly increasing function $g(s) = su'(s)$ ranging from $\alpha$ to $\beta$ on $(0,\infty)$.\\
		
		\item $c_4\neq 0$, $c_0 =0$, $c_1\neq 0$, and the polynomial $c_4x^3 +c_3 x^2 +x+c_1$ has real distict roots $\alpha,\beta,\gamma$.\\
		By Lemma \ref{lemma6.2} we have $B<\infty$ and $H(x)>0$ on $(A,B)$. It follows that all roots are real, and if we let $\alpha=A$, $\beta =B$, $\gamma<\alpha<\beta$, then we have $c_4<0$. This gives us the inequality $\gamma<-\frac{\alpha\beta}{\alpha+\beta} <0<\alpha<\beta$. 
		
		On the other hand, if we let $\alpha =A$, $\beta = B$, $\alpha<\beta<\gamma$, then we have $c_4>0$. 
		
		We can write equation \eqref{*8} as 
		\begin{align}\label{*7}
		g'(s)(\alpha\beta + \alpha\gamma +\beta\gamma) \left\{\frac{1}{(\alpha - \beta)(\alpha-\gamma)} \frac{1}{g(s)-\alpha} + 
		\frac{1}{(\beta-\alpha)(\beta-\gamma)} \frac{1}{g(s)-\beta}  \right. \\ \nonumber \left.
		+\frac{1}{(\gamma - \alpha)(\gamma-\beta)} \frac{1}{g(s)-\gamma}  
		\right\} = \frac{1}{s}.
		\end{align}
		There exists a constant $c$ such that 
		\begin{align} \tag{\ref{*10}}
		(\alpha\beta +\alpha\gamma +\beta\gamma) \left\{\frac{\log(g(s)-\alpha)}{(\alpha-\beta)(\alpha -\gamma)} + \frac{\log(\beta -g(s))}{(\beta-\alpha)(\beta -\gamma)} + \frac{\log|g(s)-\gamma|}{(\gamma-\alpha)(\gamma-\beta)} \right\}=\log s + c.
		\end{align}
		If we denote the left hand side of \eqref{*10} by $h(g(s))$, then we see that $\lim_{s\rightarrow 0^+} h(g(s)) = -\infty$, $\lim_{s\rightarrow +\infty} h(g(s)) >0$, and $\frac{d}{ds} h(g(s))>0$ on $(0,\infty)$. It follows from Lemma \ref{lemma6.1} that there exists a unique smooth strictly increasing function $g(s)$ that solves the equation.\\
		
		\item $c_4\neq 0$, $c_0\neq 0$, and the polynomial $c_4 x^4 +c_3 x^3 +x^2 +c_1 x+c_0$ has at most one real root.\\
		By Lemma \ref{lemma6.2} $B<\infty$, and the equation \eqref{ordinary differential equation} does not admit the required solution.\\
		
		\item $c_4\neq 0$, $c_0\neq 0$, and the polynomial $c_4 x^4 +c_3 x^3 +x^2 +c_1 x+c_0$ has real roots $\alpha,\alpha,\alpha,\beta$ with $\alpha<\beta$.\\
		It follows from Lemma \ref{lemma6.2} that $B<\infty$, and $H(x)>0$ on $(A,B)$. This implies $0<\alpha =A<\beta=B$, and $c_4=\frac{1}{3\alpha^2 +3\alpha\beta}<0$, which gives a contradiction.\\
		
		\item  $c_4\neq 0$, $c_0\neq 0$, and the polynomial $c_4 x^4 +c_3 x^3 +x^2 +c_1 x+c_0$ has real roots $\alpha,\beta,\beta,\beta$ with $\alpha<\beta$.\\
		It follows from Lemma \ref{lemma6.2} that $0<\alpha=A<\beta=B$, and $c_4<0$, which gives a contradiction. 
		
		\item $c_4\neq 0$, $c_0\neq 0$, and the polynomial $c_4 x^4 +c_3 x^3 +x^2 +c_1 x+c_0$ has no real roots, and has complex roots $a-ib,a-ib,a+ib,a+ib$.\\
		It follows from Lemma \ref{lemma6.2} that $H(A)=0$, which is a contradiction.\\
		
		\item $c_4\neq 0$, $c_0\neq 0$, and the polynomial $c_4 x^4 +c_3 x^3 +x^2 +c_1 x+c_0$ has real roots $\alpha,\alpha,\beta,\beta$.\\
		We have $c_4 x^4 +c_3 x^3 +x^2 +c_1 x+c_0 = c_4 (x-\alpha)^2 (x-\beta)^2$ with $\alpha\beta\neq 0$. The equation \eqref{*8} can be written as 
		\begin{align*}
		g'(s)((\alpha+\beta)^2 +2\alpha\beta)\left\{
		\frac{\alpha+\beta}{(\beta-\alpha)^3}\frac{1}{g(s)-\alpha}  +\frac{\alpha}{(\alpha-\beta)^2}\frac{1}{(g(s)-\alpha)^2} 
		\right. \\ \nonumber \left. - 
		\frac{\alpha+\beta}{(\beta-\alpha)^3}\frac{1}{g(s)-\beta} 
		+\frac{\beta}{(\alpha-\beta)^2}\frac{1}{(g(s)-\beta)^2} 
		\right\} = \frac{1}{s}.
		\end{align*} 
		We see from Lemma \ref{lemma6.2} that $\alpha = A>0$, $\beta=B$, and $c_4>0$, where $c_4 = \frac{1}{(\alpha+\beta)^2 +2\alpha\beta}$. We can integrate the above equation to obtain 
		\begin{align}\tag{\ref{*6}}
		((\alpha+\beta)^2 + 2\alpha\beta) \left\{ \frac{\alpha+\beta}{(\beta-\alpha)^3} \log(g(s)-\alpha) - \frac{\alpha}{(\alpha-\beta)^2} \frac{1}{g(s)-\alpha} - \right.\\ \nonumber \left. \frac{\alpha+\beta}{(\beta-\alpha)^3} \log(\beta - g(s)) - \frac{\beta}{(\alpha-\beta)^2} \frac{1}{g(s)-\beta} \right\} = \log s + c.
		\end{align}
		On the other hand, Lemma \ref{lemma6.1} implies that there exists a unique smooth strictly increasing function $g(s)$ ranging from $\alpha$ to $\beta$ on $(0,\infty)$.\\
		
		\item $c_4\neq 0$, $c_0\neq 0$, and the polynomial $c_4 x^4 +c_3 x^3 +x^2 +c_1 x+c_0$ has three distinct real roots $\alpha,\alpha,\beta,\gamma$ with $\alpha<\beta<\gamma$.\\
		It follows from Lemma \ref{lemma6.2} that we can either have $\alpha=A$, $\beta=B$; or $\beta=A$, $\gamma =B$.
		
		Let us start with the case $\alpha =A$, $\beta = B$. In this case we have $c_4 x^4 +c_3 x^3 +x^2 +c_1 x+c_0 = c_4(x-\alpha)^2 (x-\beta)(x-\gamma)$. Then, $c_4 =\frac{1}{\alpha^2 +2\alpha\beta +2\alpha\gamma +\beta \gamma}$ and we can see from Lemma \ref{lemma6.2} that $\alpha>0$, $c_4>0$.
		
		The equation \eqref{*8} can be rewritten as 
		\begin{align}\label{*2}
		g'(s)(\alpha^2 +2\alpha\beta+2\alpha\gamma +\beta\gamma) \left\{ \frac{\gamma}{(\gamma-\alpha)(\gamma-\beta)} \frac{1}{g(s)-\gamma} + \right. \\\nonumber \left. \frac{\alpha}{(\alpha-\gamma)(\alpha-\beta)} \frac{1}{(g(s)-\alpha)^2} + \frac{-\alpha^2+\beta\gamma}{(\alpha-\gamma)(\alpha-\beta)} \frac{1}{g(s)-\alpha} +\right. \\\nonumber \left. \frac{\beta}{(\beta-\gamma)(\beta-\alpha)^2} \frac{1}{g(s)-\beta} 
		\right\} = \frac{1}{s}.
		\end{align}
		There exists a constant $c$ such that 			
		\begin{align}\tag{\ref{*3}}
		(\alpha^2 +2\alpha\beta+2\alpha\gamma +\beta\gamma) \left\{ \frac{\gamma}{(\gamma-\alpha)(\gamma-\beta)} \log(\gamma-g(s)) - \right. \\\nonumber \left. \frac{\alpha}{(\alpha-\gamma)(\alpha-\beta)} \frac{1}{g(s)-\alpha} + \frac{-\alpha^2+\beta\gamma}{(\alpha-\gamma)(\alpha-\beta)} \log(g(s)-\alpha) +\right. \\\nonumber \left. \frac{\beta}{(\beta-\gamma)(\beta-\alpha)^2} \log|\beta-g(s)| 
		\right\} = \log s + c.
		\end{align}
		Note that $-\alpha^2+\beta\gamma>0$. By Lemma \ref{lemma6.1}, there exists a unique smooth function $g(s):(0,\infty)\rightarrow (\alpha,\beta)$ with $g'(s)>0$ which solves the above equation.
		
		On the other hand, if we assume $\beta=A$ and $\gamma=B$, then it follows from Lemma \ref{lemma6.2} that $c_4 = \frac{1}{\alpha^2 +2\alpha \beta +2\alpha\gamma+\beta\gamma}<0$, and $0<\beta<\gamma$. 
		
		Equivalently, we can write $\alpha<-\frac{\alpha^2+\beta\gamma}{2(\beta+\gamma)}<0<\beta<\gamma$. Note that for any given $0<\beta<\gamma$, such $\alpha$ values exist.
		
		Equation \eqref{*8} can be rewritten as equation \eqref{*2} as before, however, this time we are looking for a smooth soution $g(s)$ with values in $(\beta,\gamma)$. Keeping this in mind, we investigate \eqref{*2} to obtain \eqref{*3}, and use Lemma \ref{lemma6.1} to conclude that there exists a unique smooth strictly increasing function $g:(0,\infty)\rightarrow(\beta,\gamma)$ satisfying \eqref{*3}.\\
		
		\item $c_4\neq 0$, $c_0\neq 0$, and the polynomial $c_4 x^4 +c_3 x^3 +x^2 +c_1 x+c_0$ has three distinct real roots $\alpha,\beta,\beta,\gamma$ with $\alpha<\beta<\gamma$.\\
		The equation \eqref{ordinary differential equation} can be rewritten as
		\begin{align*}
		g'(s)(\beta^2 +2\beta\alpha+2\beta\gamma +\alpha\gamma)  \left\{ \frac{\alpha}{(\alpha-\beta)(\alpha -\gamma)} \frac{1}{g(s)-\alpha} \right.\\ \nonumber \left. + \frac{\beta}{(\beta-\alpha)(\beta-\gamma)} \frac{1}{(g(s)-\beta)^2} + \frac{-\beta^2+\alpha\gamma}{(\beta-\alpha)(\beta-\gamma)} \frac{1}{g(s)-\beta}  \right. \\\nonumber \left. + \frac{\gamma}{(\gamma-\alpha)(\gamma-\beta)^2} \frac{1}{g(s)-\gamma} 
		\right\} = \frac{1}{s}.
		\end{align*}
		There exists a constant $c$ such that
		\begin{align}\tag{\ref{*4}}
		(\beta^2 +2\beta\alpha+2\beta\gamma +\alpha\gamma) \left\{ \frac{\alpha}{(\alpha-\beta)(\alpha -\gamma)} \log(g(s)-\alpha) \right.\\ \nonumber \left. - \frac{\beta}{(\beta-\alpha)(\beta-\gamma)} \frac{1}{g(s)-\beta} + \frac{-\beta^2+\alpha\gamma}{(\beta-\alpha)(\beta-\gamma)} \log|g(s)-\beta| \right. \\\nonumber \left. + \frac{\gamma}{(\gamma-\alpha)(\gamma-\beta)^2} \log(\gamma-g(s)) 
		\right\} = \log s + c.
		\end{align}
		It follows from Lemma \ref{lemma6.2} that we have $B<\infty$ for degree reasons, so we can choose either $\alpha=A$, $\beta=B$; or $\beta=A$, $\gamma=B$. By Lemma \ref{lemma6.2} we have $H(x)>0$ on $(A,B)$, which implies $c_4<0$ in both cases. However, since $c_4=\frac{1}{\beta^2+2\beta\alpha+2\beta\gamma+\alpha\gamma}$ and $A>0$,  we see that the former case is impossible, leaving us with the choice $\beta=A$, $\gamma = B$. It follows from Lemma \ref{lemma6.1} that there exists a unique smooth strictly increasing function $g:(0,\infty)\rightarrow(\beta,\gamma)$ satisfying \eqref{*4}. \\
		
		\item  $c_4\neq 0$, $c_0\neq 0$, and the polynomial $c_4 x^4 +c_3 x^3 +x^2 +c_1 x+c_0$ has three distinct real roots $\alpha,\beta,\gamma,\gamma$ with $\alpha<\beta<\gamma$.\\
		It follows from Lemma \ref{lemma6.2} that we can have either $\alpha=A$, $\beta=B$; or $\beta=A$, $\gamma=B$. In the former case, Lemma \ref{lemma6.2} implies $c_4=\frac{1}{\alpha\beta+2\alpha\gamma + 2\beta\gamma +\gamma^2}<0$, which contradicts with our choice $0<\alpha=A<\beta=B<\gamma$.
		
		Let us assume $\beta=A$, $\gamma=B$. Since $H(x)>0$ on $(\beta,\gamma)$, we have $c_4>0$, which implies that $-\frac{\gamma(\gamma+2\beta)}{2\gamma+\beta}<\alpha$. We have $\alpha\neq 0$ as $c_0\neq 0$. The equation \eqref{*8} can be written as
		\begin{align*}
		g'(\gamma^2 +2\beta\gamma+2\alpha\gamma +\alpha\beta) \left\{ \frac{\alpha}{(\alpha-\beta)(\alpha -\gamma)} \frac{1}{g(s)-\alpha} \right.\\ \nonumber \left. + \frac{\gamma}{(\gamma-\alpha)(\gamma-\beta)} \frac{1}{(g(s)-\gamma)^2} + \frac{-\gamma^2+\alpha\beta}{(\gamma-\alpha)(\gamma-\beta)} \frac{1}{g(s)-\gamma} \right. \\\nonumber \left. + \frac{\beta}{(\beta-\alpha)(\beta-\gamma)^2} \frac{1}{g(s)-\beta}
		\right\} = \frac{1}{s}.
		\end{align*}
		There exists a constant $c$ such that
		\begin{align}\tag{\ref{*5}}
		(\gamma^2 +2\beta\gamma+2\alpha\gamma +\alpha\beta) \left\{ \frac{\alpha}{(\alpha-\beta)(\alpha -\gamma)} \log(g(s)-\alpha) \right.\\ \nonumber \left. - \frac{\gamma}{(\gamma-\alpha)(\gamma-\beta)} \frac{1}{g(s)-\gamma} + \frac{-\gamma^2+\alpha\beta}{(\gamma-\alpha)(\gamma-\beta)} \log(\gamma-g(s)) \right. \\\nonumber \left. + \frac{\beta}{(\beta-\alpha)(\beta-\gamma)^2} \log(g(s)-\beta) 
		\right\} = \log s + c.
		\end{align}
		It follows from Lemma \ref{lemma6.1} that there exists a unique smooth strictly increasing function $g(s):(0,\infty)\rightarrow(\beta,\gamma)$ that solves equation \eqref{*5}.\\
		
		\item  $c_4\neq 0$, $c_0\neq 0$, and the polynomial $c_4 x^4 +c_3 x^3 +x^2 +c_1 x+c_0$ has four distinct real roots $\alpha,\beta,\gamma,\tau$ with $\alpha<\beta<\gamma<\tau$, and $\alpha\beta\gamma\tau\neq 0$.\\
		Equation \eqref{*8} can be rewritten as 
		\begin{align*}
		g'(\alpha\beta+\alpha\gamma+\alpha\tau+\beta\gamma+\beta\tau+\gamma\tau) \left\{ \frac{\alpha}{(\alpha-\beta)(\alpha-\gamma)(\alpha-\tau)}\frac{1}{g(s)-\alpha}  \right. \\\nonumber \left. +
		\frac{\beta}{(\beta-\alpha)(\beta-\gamma)(\beta-\tau)}\frac{1}{g(s)-\beta} +
		\frac{\gamma}{(\gamma-\alpha)(\gamma-\beta)(\gamma-\tau)}\frac{1}{g(s)-\gamma} \right. \\\nonumber \left. +
		\frac{\tau}{(\tau-\alpha)(\tau-\beta)(\tau-\gamma)}\frac{1}{g(s)-\tau} 
		\right\}
		=\frac{1}{s} .
		\end{align*}
		There exists a constant $c$ such that 
		\begin{align}\tag{\ref{**6}}
		(\alpha\beta+\alpha\gamma+\alpha\tau+\beta\gamma+\beta\tau+\gamma\tau) \left\{ \frac{\alpha}{(\alpha-\beta)(\alpha-\gamma)(\alpha-\tau)}\log(g(s)-\alpha)  \right. \\\nonumber \left. +
		\frac{\beta}{(\beta-\alpha)(\beta-\gamma)(\beta-\tau)}\log|g(s)-\beta|+
		\frac{\gamma}{(\gamma-\alpha)(\gamma-\beta)(\gamma-\tau)}\log|g(s)-\gamma| \right. \\\nonumber \left. +
		\frac{\tau}{(\tau-\alpha)(\tau-\beta)(\tau-\gamma)}\log|g(s)-\tau| 
		\right\}
		=\log s + c .
		\end{align}
		It follows from Lemma \ref{lemma6.2} that we have $A>0$, $B<\infty$, and $H(x)>0$ on $(A,B)$. This implies that we have three possibilities:
		\begin{enumerate}[label=(\roman*)]
			\item $\alpha=A$, $\beta=B$ and $c_4<0$.\\
			In this case, all roots of the polynomial $H(x)$ are positive which contradicts with $c_4<0$.
			\item $\beta=A$, $\gamma=B$, and $c_4>0$\\
			By Lemma \ref{lemma6.1}, there exists a unique smooth, strictly increasing function $g(s):(0,\infty)\rightarrow(\beta,\gamma)$ satisfying equation \eqref{**6}, whenever $\alpha>-\frac{\beta\gamma+\beta\tau+\gamma\tau}{\beta+\gamma+\tau}$.
			\item $\gamma=A$, $\beta=\tau$, and $c_4<0$.\\
			By Lemma \ref{lemma6.1}, there exists a unique smooth, strictly increasing function $g(s):(0,\infty)\rightarrow(\gamma,\tau)$ satisfying equation \eqref{**6}, whenever $\alpha\beta+\alpha\gamma+\alpha\tau+\beta\gamma+\beta\tau+\gamma\tau<0$.\\
		\end{enumerate}

		\item  $c_4\neq 0$, $c_0\neq 0$, and the polynomial $c_4 x^4 +c_3 x^3 +x^2 +c_1 x+c_0$ has four distinct  roots $\alpha,\beta,a+ib,a-ib$.\\
		It follows from Lemma \ref{lemma6.2} that $\alpha=A$, $\beta=B$, and $c_4<0$. If we write $H(x)=c_4 (x-\alpha)(x-\beta)(x^2+2ax+a^2+b^2)$, then $c_4<0$ can be written as $a^2+2a(\alpha+\beta) + b^2+\alpha\beta<0$. This condition holds for those $\alpha,\beta,a,b$ which satisfy $b^2<\alpha^2+\alpha\beta+\beta^2$ and $-(\alpha+\beta)-\sqrt{\alpha^2 +\alpha\beta+\beta^2-b^2}<a<-(\alpha+\beta) + \sqrt{\alpha^2+\alpha\beta+\beta^2-b^2}$.
		
		The equation \eqref{*8} can be rewritten as
		\begin{align*}
		(a^2 +2a(\alpha+\beta)+b^2+\alpha\beta)g' \left\{c_1 \frac{c_1}{g(s)-\alpha}  +\frac{c_2}{g(s)-\beta}  +  \right. \\\nonumber \left. 
		\frac{-(c_1 + c_2) g(s) + 2(c_1+c_2) a - c_1 \alpha -c_2 \beta}{g^2(s) -2a g(s)+a^2+b^2} \right\} = \frac{1}{s},
		\end{align*}
		where $c_1 = \frac{\alpha}{(\alpha-\beta)(\alpha^2-2a\alpha +a^2+b^2)}$ and $c_2 = \frac{\beta}{(\beta -\alpha)(\beta^2-2a\beta +a^2+b^2)}$.
		
		On the other hand, Lemma \ref{lemma6.1} implies that there exists a unique smooth, strictly increasing function $g(s):(0,\infty)\rightarrow (\alpha,\beta)$ determined by 
		\begin{align}\tag{\ref{**7}}
		(a^2 +2a(\alpha+\beta)+b^2+\alpha\beta) \left\{c_1 \log(g(s)-\alpha) +c_2 \log(\beta -g(s)) +  \right. \\\nonumber \left. \int_{\alpha}^{\beta}\frac{-(c_1 + c_2) g(s) + 2(c_1+c_2) a - c_1 \alpha -c_2 \beta}{g^2(s) -2a g(s)+a^2+b^2} g'(s) ds\right\} = \log s.
		\end{align}
		Here we note that the integral in equation \eqref{**7} is a proper integral, since the denominator is never zero.
	\end{enumerate}
\end{proof}

\section{Technical Lemmas}\label{section technical}

In this section we include the technical lemmas by other authors that we use to prove our results.

\begin{proposition}[Monn \cite{Monn85}, Proposition 2.1]\label{MonnProp2.1}
	Let $B$ be an open ball containing the origin in $\mathbb{C}^n$. Let $u$ be a radial function on $\overline B$, and let $\tilde u(r)=u(r,0,\dots,0)$. Then $u\in C^k(\overline B)$ if and only if $\tilde u\in C^k[0,1]$, and $\tilde u^{(\ell)}(0)=0$ for all $\ell\leq k$, $\ell$ odd.
\end{proposition}

\begin{proposition}[Monn \cite{Monn85}, Proposition 4.1]\label{MonnProp4.1}
	The $k^{\normalfont th}$ derivative of two real-valued functions, $f\circ g$, can be written as a sum of terms of the form 
	$$f^{(\lambda)} (g)\cdot P(g',g'',\dots,g^{k+1-\lambda})$$ 
	where $P$ is a monomial of degree $\lambda\leq k$ and of weighted degree $k$.
\end{proposition}

The following lemma is useful for eliminating impossible cases as solutions of the extremal equation $sg'=F(g)$.

\begin{lemma}[\cite{HeLi}, Lemma 6.2]\label{lemma6.2} Suppose $H(x)$ is a polynomial of degree $m$ and the ordinary differential equation $$ \frac{g^k(s)g'(s)}{H(g(s))} =\frac{1}{s}$$ admits a smooth solution $g(s)$ on $(0,\infty)$ with $g(s)>0$, $g'(s)>0$. Denote by $A=\displaystyle\lim_{s\rightarrow 0^+} g(s)$, $B=\displaystyle\lim_{s\rightarrow +\infty} g(s)$. Then
	\begin{enumerate}
		\item $H(A)=0$ and $H(x)>0$ for $x\in(A,B)$.
		\item If $B=+\infty$, then $\deg H\leq k+1$. Moreover, $A$ is the largest and nonnegative real root of $H(x)$, and $H(x)$ is positive on $(A,+\infty)$.
		\item If $B < +\infty$, then $H(B)=0$. Moreover, $A$ and $B$ are two successive nonnegative real roots of the polynomial $H(x)$, and $H(x)$ is positive on the interval $(A,B)$
	\end{enumerate}
\end{lemma}

Once the impossible cases are eliminated by the above lemma, we use the following lemma to show the existence of solutions.
\begin{lemma}[\cite{HeLi}, Lemma 6.1]\label{lemma6.1}
	Let $h:(A,B)\rightarrow \mathbb{R}$ be a smooth, strictly increasing function with $\displaystyle\lim_{t\rightarrow A} h(t)= -\infty$, $\displaystyle\lim_{t\rightarrow B} h(t)= \infty$. Then, for any constant $a>0$ and $c$, there exists a unique smooth, strictly increasing function $g:(0,\infty)\rightarrow\mathbb{R}$ such that $$h(g(s))=a\log s + c$$
	and $\displaystyle\lim_{s\rightarrow 0^+} g(s)= A$, $\displaystyle\lim_{s\rightarrow +\infty} g(s)= B$.
\end{lemma}

\bibliographystyle{alpha}
\bibliography{taskent_arxiv_1}

\begin{thebibliography}{Mon86}

\bibitem[Abr01]{Abreu01}
Miguel Abreu.
\newblock K\"{a}hler metrics on toric orbifolds.
\newblock {\em J. Differential Geom.}, 58(1):151--187, 2001.

\bibitem[Abr10]{Abreu10}
Miguel Abreu.
\newblock Toric {K}\"{a}hler metrics: cohomogeneity one examples of constant
  scalar curvature in action-angle coordinates.
\newblock {\em J. Geom. Symmetry Phys.}, 17:1--33, 2010.

\bibitem[AG02]{AG02}
Vestislav Apostolov and Paul Gauduchon.
\newblock Selfdual {E}instein {H}ermitian four-manifolds.
\newblock {\em Ann. Sc. Norm. Super. Pisa Cl. Sci. (5)}, 1(1):203--243, 2002.

\bibitem[AR17]{AR17}
Vestislav Apostolov and Yann Rollin.
\newblock A{LE} scalar-flat {K}\"{a}hler metrics on non-compact weighted
  projective spaces.
\newblock {\em Math. Ann.}, 367(3-4):1685--1726, 2017.

\bibitem[Bry01]{B01}
Robert~L. Bryant.
\newblock Bochner-{K}\"{a}hler metrics.
\newblock {\em J. Amer. Math. Soc.}, 14(3):623--715, 2001.

\bibitem[Cal82]{Calabi82}
Eugenio Calabi.
\newblock Extremal {K}\"{a}hler metrics.
\newblock In {\em Seminar on {D}ifferential {G}eometry}, volume 102 of {\em
  Ann. of Math. Stud.}, pages 259--290. Princeton Univ. Press, Princeton, N.J.,
  1982.

\bibitem[Cao96]{C94}
Huai-Dong Cao.
\newblock Existence of gradient {K}\"{a}hler-{R}icci solitons.
\newblock In {\em Elliptic and parabolic methods in geometry ({M}inneapolis,
  {MN}, 1994)}, pages 1--16. A K Peters, Wellesley, MA, 1996.

\bibitem[Cao97]{C97}
Huai-Dong Cao.
\newblock Limits of solutions to the {K}\"{a}hler-{R}icci flow.
\newblock {\em J. Differential Geom.}, 45(2):257--272, 1997.

\bibitem[DL16]{DL16}
Michael~G. Dabkowski and Michael~T. Lock.
\newblock On {K}\"{a}hler conformal compactifications of {$U(n)$}-invariant
  {ALE} spaces.
\newblock {\em Ann. Global Anal. Geom.}, 49(1):73--85, 2016.

\bibitem[EH79]{EH79}
Tohru Eguchi and Andrew~J. Hanson.
\newblock Self-dual solutions to {E}uclidean gravity.
\newblock {\em Ann. Physics}, 120(1):82--106, 1979.

\bibitem[FIK03]{FIK03}
Mikhail Feldman, Tom Ilmanen, and Dan Knopf.
\newblock Rotationally symmetric shrinking and expanding gradient
  {K}\"{a}hler-{R}icci solitons.
\newblock {\em J. Differential Geom.}, 65(2):169--209, 2003.

\bibitem[Gau19]{G19}
Paul Gauduchon.
\newblock Invariant scalar-flat {K}\"{a}hler metrics on {$\mathrm{Cal}\,
  O(-\ell)$}.
\newblock {\em Boll. Unione Mat. Ital.}, 12(1-2):265--292, 2019.

\bibitem[HL18]{HeLi}
Weiyong {He} and Jun {Li}.
\newblock {Rotation invariant singular K\"ahler metrics with constant scalar
  curvature on $\mathbb{C}^n$}.
\newblock {\em arXiv e-prints}, page arXiv:1808.03925, August 2018.

\bibitem[HS97]{HS97}
Andrew~D. Hwang and Santiago~R. Simanca.
\newblock Extremal {K}\"{a}hler metrics on {H}irzebruch surfaces which are
  locally conformally equivalent to {E}instein metrics.
\newblock {\em Math. Ann.}, 309(1):97--106, 1997.

\bibitem[HS02]{HS02}
Andrew~D. Hwang and Michael~A. Singer.
\newblock A momentum construction for circle-invariant {K}\"{a}hler metrics.
\newblock {\em Trans. Amer. Math. Soc.}, 354(6):2285--2325, 2002.

\bibitem[Kle77]{K77}
Paul~F. Klembeck.
\newblock A complete {K}\"{a}hler metric of positive curvature on {$C^{n}$}.
\newblock {\em Proc. Amer. Math. Soc.}, 64(2):313--316, 1977.

\bibitem[LeB88]{LeBrun88}
Claude LeBrun.
\newblock Counter-examples to the generalized positive action conjecture.
\newblock {\em Comm. Math. Phys.}, 118(4):591--596, 1988.

\bibitem[Mon86]{Monn85}
David Monn.
\newblock Regularity of the complex {M}onge-{A}mpere equation for radially
  symmetric functions of the unit ball.
\newblock {\em Math. Ann.}, 275(3):501--511, 1986.

\bibitem[PP91]{PP91}
Henrik Pedersen and Y.~Sun Poon.
\newblock Hamiltonian constructions of {K}\"{a}hler-{E}instein metrics and
  {K}\"{a}hler metrics of constant scalar curvature.
\newblock {\em Comm. Math. Phys.}, 136(2):309--326, 1991.

\bibitem[Sim91]{S91}
Santiago~R. Simanca.
\newblock K\"{a}hler metrics of constant scalar curvature on bundles over
  {${\bf C}{\rm P}_{n-1}$}.
\newblock {\em Math. Ann.}, 291(2):239--246, 1991.

\bibitem[TL70]{TL70}
Shun-ichi Tachibana and Richard~Chieng Liu.
\newblock Notes on {K}\"{a}hlerian metrics with vanishing {B}ochner curvature
  tensor.
\newblock {\em K{o}dai Math. Sem. Rep.}, 22:313--321, 1970.

\bibitem[WZ11]{WZ11}
Hung-Hsi Wu and Fangyang Zheng.
\newblock Examples of positively curved complete {K}\"{a}hler manifolds.
\newblock In {\em Geometry and analysis. {N}o. 1}, volume~17 of {\em Adv. Lect.
  Math. (ALM)}, pages 517--542. Int. Press, Somerville, MA, 2011.

\end{thebibliography}

\vspace{40pt}

\textit{email: selintaskent@gmail.com}

\end{document}